\documentclass[aap,MSNbibl,nameyear,dvips]{arximspdf}
\usepackage{dcolumn}
\usepackage{graphicx}

\doi{10.1214/11-AAP819} 
\volume{22}
\issue{5}
\pubyear{2012}
\firstpage{1962}
\lastpage{1988}

\makeatletter
\newcolumntype{d}[1]{D{.}{.}{#1}}
\newcommand{\eqref}[1]{(\ref{#1})}
\newcommand{\A}{A}
\newcommand{\Ai}{A}

\newtheorem{lemma}{Lemma}
\newtheorem{theorem}{Theorem}
\newtheorem{proposition}{Proposition}
\newtheorem{corollary}{Corollary}
\makeatother

\begin{document}
\begin{frontmatter}

\title{Fast approach to the Tracy--Widom law at the~edge~of~GOE~and~GUE\thanksref{T1}}
\runtitle{Convergence rate for GUE and GOE}
\thankstext{T1}{Supported in part by DMS-09-06812 and NIH BIB R01EB1988.}

\begin{aug}
\author[A]{\fnms{Iain M.} \snm{Johnstone}\corref{}\ead[label=e1]{imj@stanford.edu}}
\and
\author[B]{\fnms{Zongming} \snm{Ma}\ead[label=e2]{zongming@wharton.upenn.edu}}
\runauthor{I. M. Johnstone and Z. Ma}
\affiliation{Stanford University and University of Pennsylvania}
\address[A]{Department of Statistics\\
Stanford University\\
Sequoia Hall\\
Stanford, California 94305\\
USA\\
\printead{e1}} 
\address[B]{Department of Statistics\\
The Wharton School\\
University of Pennsylvania\\
Philadelphia, Pennsylvania 19104\\
USA\\
\printead{e2}}
\end{aug}

\received{\smonth{8} \syear{2010}}
\revised{\smonth{9} \syear{2011}}

%
\begin{abstract}
We study the rate of convergence for the largest eigenvalue
distributions in the Gaussian unitary and orthogonal ensembles to their
Tracy--Widom limits.

We show that one can achieve an $O(N^{-2/3})$ rate with particular
choices of the centering and scaling constants. The arguments here also
shed light on more complicated cases of Laguerre and Jacobi ensembles,
in both unitary and orthogonal versions.

Numerical work shows that the suggested constants yield reasonable
approximations, even for surprisingly small values of $N$.
\end{abstract}

%

\end{frontmatter}

\section{Introduction}
\label{secintroduction}

The celebrated papers of \citeauthor{trwi94} (\citeyear{trwi94,trwi96})
described the limiting
distributions of the largest eigenvalues of the Gaussian unitary and
orthogonal ensembles (GUE and GOE), respectively.
The purpose of this article is to show that an appropriate choice of
centering and scaling allows us to establish a rate of convergence
in these results, and further that this rate can be understood as
``second order,'' being $O(N^{-2/3})$ rather than the $O(N^{-1/3})$
that would otherwise apply.

The Gaussian ensembles refer, as is usual, to eigenvalue densities of
$x = (x_1, \ldots, x_N)$ given by
%
\[
f(x) = c_{N\beta} \prod_{i=1}^N e^{- \beta x_i^2/2} \prod_{i<j} |x_i
- x_j|^\beta,
\]
with $\beta= 1$ corresponding to $\operatorname{GOE}_N$ and $\beta=
2$ to
$\operatorname{GUE}_N$, the subscript being shown only when clarity dictates.
The corresponding matrix models specify that~$f$ is the density of the
eigenvalues $x$ of a symmetric or Hermitian random matrix~$M$ with
independent entries on and above the diagonal, whose density function is
given by
%
\[
g(M) = c_{N \beta}' \exp\{ - (\beta/ 2) \operatorname{tr} M^2 \}.
\]

Our principal rate of convergence results follow. The Tracy--Widom
distributions are denoted $F_\beta(s)$ for $\beta= 1, 2$.

\begin{theorem} \label{thunitary-main-res}
Let $x_{(1)}$ denote the largest eigenvalue of a sample from
$\operatorname{GUE}_N$,
and
%
\[
\mu_N = \sqrt{2N},\qquad
\tau_N = 2^{-1/2} N^{-1/6}.
\]
Given $s_0$, there exists $C = C(s_0)$
such that for $s \geq s_0,$
%
\[
\bigl|P\bigl\{ \bigl( x_{(1)} - \mu_N\bigr)/\tau_N \leq s \bigr
\} - F_2(s)\bigr|
\leq C N^{-2/3} e^{-s}.
\]
\end{theorem}




\begin{theorem} \label{thorthogonal-main-res}
Let $x_{(1)}$ denote the largest eigenvalue of a sample from
$\operatorname{GOE}_{N+1}$,
with $N+1$ even, and
%
\begin{equation}
\label{eqrealcasemusig}
\mu_N = \sqrt{2N +1},\qquad\tau_N = 2^{-1/2} N^{-1/6}.
\end{equation}
Given $s_0$, there exists $C = C(s_0)$
such that for $s \geq s_0,$
%
\begin{equation}
\label{eqgoe-main-bd}
\bigl|P\bigl\{ \bigl( x_{(1)} - \mu_N\bigr)/\tau_N \leq s \bigr
\} - F_1(s)\bigr|
\leq C N^{-2/3} e^{-s/2}.
\end{equation}
\end{theorem}

We use index $N+1$ (rather than $N$) because of a key formula
relating the Gaussian orthogonal ensemble $\operatorname{GOE}_{N+1}$
to the Gaussian unitary
ensemble $\operatorname{GUE}_N$, (\ref{eqRealKernel0}) below. The
centering and
scaling constants carry subscripts $N$ rather than $N+1$ for this reason.

Our interest in these results is threefold.
First, they provide the simplest case of a class of such $O(N^{-2/3})$
convergence results for the classical orthogonal polynomial
ensembles---the other two being the
Laguerre and Jacobi ensembles---in both orthogonal and
unitary versions.
These results are of interest in statistics because they show that the
Tracy--Widom approximation is accurate enough to replace exact
evaluation of the finite LOE and JOE probabilities for many applied
purposes where highly accurate values are not necessary
[\citet{john08a}].
The results of this paper focus on the corresponding phenomenon for
the simplest case of GUE and GOE.
Since the LOE and JOE proofs are lengthy analyses with Laguerre and
Jacobi polynomial asymptotics, respectively, this paper outlines
the approach in the simplest case.

Second, our interest was stimulated by \citet{chou09}, which provided
the leading terms in an Edgeworth expansion of the largest eigenvalue
distribution of GOE, and remarked that the $N^{-1/3}$
correction term does not
vanish in GOE.
As our earlier results on $O(N^{-2/3})$ convergence for LOE and JOE
would suggest, a similar $O(N^{-2/3})$ property for GOE with a suitable
specific centering, it seemed,
therefore, of interest to verify the conjecture in this setting.
Although we subsequently learned of an error in the argument of
\citet{chou09} (private communication), it was an important stimulus
for this work.



Third, we find it of interest that adjustment of $\mu_N$ and $\tau_N$
to secure $O(N^{-2/3})$ convergence yields an approximation, that is,
adequate---for some purposes---for surprisingly small values of $N$.

To illustrate, first in GUE, Table~\ref{tabgue} shows the exact
probabilities $P\{ x_{(1)} \leq\mu_N + \tau_N s_{2 \alpha} \}$ for
quantiles $s_{2 \alpha}$ of the limiting $F_2$ distribution, computed
using the finite GUE function
provided in the \textsc{Matlab} toolbox \texttt{RMTFredholm}
[\citet{born09}].

%
\begin{table}
\caption{GUE approximation. For each percentile $\alpha$ shown in
the top row, let $s_{2 \alpha}$ be the quantile $F_2(s_{2
\alpha}) = \alpha$, and $s_{N \alpha} = \mu_N + \tau_N
s_{2 \alpha}$ for $\mu_N$ and $\tau_N$ specified in Theorem
\protect\ref{thunitary-main-res}. Table entries are $P \{ x_{(1)}
\leq
s_{N \alpha} \}$ computed using Bornemann's code for $E_2^{(n)}(0;
[s_{N \alpha}, \infty) )$}\label{tabgue}
\begin{tabular*}{\textwidth}{@{\extracolsep{\fill}}ld{2.3}ccccccccc@{}}
\hline
$\bolds{N}$&\multicolumn{1}{c}{$\bolds{\mu_N}$}&\textbf
{0.01}&\textbf{0.05}&\textbf{0.1}&\textbf{0.3}&\textbf{0.5}&\textbf
{0.7}&\textbf{0.9}&\textbf{0.95}&\textbf{0.99}\\
\hline
\phantom
{00}2&2.000&0.026&0.087&0.149&0.359&0.549&0.732&0.913&0.958&0.992\\
\phantom
{00}5&3.162&0.017&0.068&0.125&0.331&0.526&0.717&0.907&0.954&0.991\\
\phantom
{0}10&4.472&0.014&0.061&0.115&0.319&0.516&0.711&0.905&0.953&0.991\\
\phantom
{0}25&7.071&0.012&0.056&0.108&0.310&0.509&0.706&0.902&0.951&0.991\\
\phantom
{0}50&10.000&0.011&0.054&0.105&0.307&0.506&0.704&0.902&0.951&0.990\\
\phantom
{0}75&12.247&0.011&0.053&0.104&0.305&0.504&0.703&0.901&0.951&0.990\\
100&14.142&0.011&0.052&0.103&0.304&0.504&0.702&0.901&0.951&0.990\\
200&20.000&0.011&0.051&0.102&0.303&0.502&0.701&0.901&0.950&0.990\\
500&31.623&0.010&0.051&0.101&0.301&0.501&0.701&0.900&0.950&0.990\\
\hline
\end{tabular*}
%
\end{table}

In fact, our proof suggests a slightly different centering value,
$\mu_N =  (\sqrt{2N - 1} + \sqrt{2N + 1})/2$, which differs from
$\sqrt{2N}$ in relative terms by only $O(N^{-4})$.
However, Table~\ref{tabguemod} shows an observable improvement at very
small values of $N$.

%
\begin{table}[b]
\caption{GUE approximation: as for Table \protect\ref{tabgue}, but
with $\mu_N = (\sqrt{2N-1} + \sqrt{2N+1})/2$}
\label{tabguemod}
\begin{tabular*}{\textwidth}{@{\extracolsep{\fill}}ld{2.3}ccccccccc@{}}
\hline
$\bolds{N}$&\multicolumn{1}{c}{$\bolds{\mu_N}$}&\textbf
{0.01}&\textbf{0.05}&\textbf{0.1}&\textbf{0.3}&
\textbf{0.5}&\textbf{0.7}&\textbf{0.9}&\textbf{0.95}&\textbf
{0.99}\\
\hline
\phantom
{0}2&1.984&0.025&0.083&0.143&0.349&0.538&0.723&0.909&0.955&0.992\\
\phantom
{0}5&3.158&0.017&0.067&0.123&0.328&0.523&0.715&0.906&0.954&0.991\\
10&4.471&0.014&0.061&0.115&0.318&0.515&0.710&0.904&0.952&0.991\\
25&7.071&0.012&0.056&0.108&0.310&0.508&0.706&0.902&0.951&0.990\\
50&10.000&0.011&0.054&0.105&0.306&0.505&0.704&0.901&0.951&0.990\\
\hline
\end{tabular*}
%
%
\end{table}

Our interest is primarily with GOE, for which software for
exact computation appears to be as yet unavailable.
Table~\ref{tabgoe} shows Monte Carlo simulations of
$P\{ x_{(1)} \leq\mu_N + \tau_N s_{1 \alpha} \}$ for quantiles $s_{1
\alpha}$
of the $F_1$ limit, based on $R = 10^6$ replications.
The corresponding $95 \%$ confidence intervals have half-width
$2 \sqrt{p_\alpha(1- p_\alpha)} \times10^{-3}$ which decreases from
$0.001$ at $p_\alpha= 0.5$ to $0.0002$ at $p_\alpha= 0.01$ and
$0.99$.
Thus the tabulated values should be correct to within $\pm0.001.$

%
\begin{table}
\caption{GOE approximation. Let quantiles $s_{1 \alpha}$ be defined by
$F_1(s_{1 \alpha}) = \alpha$ for values of $\alpha$ shown in the top
row, and $\mu_N$
and $\tau_N$ be given by (\protect\ref{eqrealcasemusig}).
Based on $R = 10^6$ replications drawn from GOE,
table entries are the fraction of replications of
$s_{(1)} = (x_{(1)} - \mu_N)/\tau_N$ satisfying
$s_{(1)} \leq s_{1 \alpha}$}
\label{tabgoe}
\begin{tabular*}{\textwidth}{@{\extracolsep{\fill}}lccccccccc@{}}
\hline
$\bolds{N+1}$&\textbf{0.01}&\textbf{0.05}&\textbf{0.1}&\textbf
{0.3}&\textbf{0.5}&\textbf{0.7}&\textbf{0.9}&\textbf{0.95}&\textbf
{0.99}\\
\hline
\phantom{00}2&0.010&0.045&0.090&0.279&0.483&0.698&0.914&0.963&0.995\\
\phantom{00}5&0.012&0.053&0.103&0.300&0.500&0.704&0.907&0.956&0.993\\
\phantom{0}10&0.011&0.053&0.103&0.302&0.502&0.703&0.904&0.954&0.992\\
\phantom{0}25&0.011&0.052&0.103&0.302&0.502&0.702&0.902&0.952&0.991\\
\phantom{0}50&0.011&0.052&0.102&0.301&0.501&0.701&0.901&0.951&0.991\\
\phantom{0}75&0.011&0.051&0.102&0.302&0.501&0.701&0.901&0.951&0.990\\
100&0.010&0.051&0.101&0.301&0.501&0.701&0.901&0.951&0.990\\
200&0.011&0.051&0.101&0.301&0.501&0.701&0.901&0.951&0.990\\
500&0.010&0.050&0.100&0.300&0.500&0.700&0.901&0.951&0.990\\
\hline
\end{tabular*}
\end{table}


Two features of the numerical results deserve note.
First the approximations are somewhat better in the near right tail
than in the left.
This is presumably because the underlying approximation of Hermite
polynomials by the Airy function is anchored at the turning point $0$
of the Airy equation $A^{\prime\prime}(s) = s A(s)$, which lies
in the right tail at about the $83$rd percentile of $F_1$ and the
$97$th percentile of $F_2$.

\begin{figure}[b]

\includegraphics{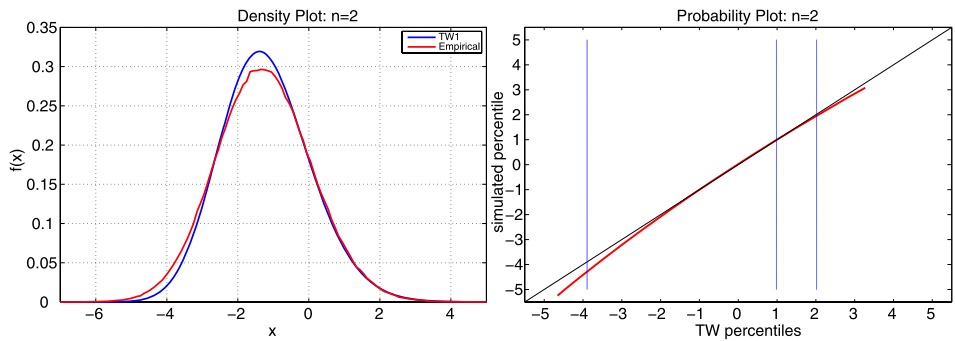}

\caption{Left panel: estimated density function for $s_{(1)}^\prime=
(x_{(1)} - \mu_N^\prime(1/5))/\tau_N(1)$, compare~(\protect\ref{eqmodcentsc}),
based on $R = 10^6$ samples from
$\operatorname{GOE}_2$,
compared to the Tracy--Widom density function~$f_1$.
Right panel: probability plot of percentiles of the $F_1$ distribution
on horizontal axis versus ordered values of $s_{(1)}$ on the vertical
axis, based on the same
samples from $\operatorname{GOE}_2$. Vertical lines mark the 1st, 95th and
99th percentiles of $F_1$.}\label{figgraphs}
\end{figure}

Second, the errors in Tables~\ref{tabgue}--\ref{tabgoe} all have the
same sign, suggesting that a~further shift in the approximating
distribution might improve accuracy. We experimented in GOE with small
changes of the form, setting $N_+ = N + 1/2$,
%
\begin{equation}
\label{eqmodcentsc}
\mu_N(\gamma) = (2N_+ - \gamma N_+^{-1/3})^{1/2},
\qquad
\tau_N(c) = 2^{-1/2} (N+c)^{-1/6},
\end{equation}
and obtained good results, Table~\ref{tabgoemod} and Figure
\ref{figgraphs}, for $\gamma= 1/5$ and $c = 1$.

These values differ from $\mu_N$ and $\tau_N$ of Theorem
\ref{thorthogonal-main-res} by
by relative errors of $O(N^{-4/3})$ and $O(N^{-1})$, respectively,
and so have no
effect on the validity of Theorem~\ref{thorthogonal-main-res}.
However, they provide a substantial numerical improvement, especially in
the right tail for values of $N$ below $10.$
Indeed, for some purposes, the approximation in the right tail would
be adequate, even for $N = 2$.

%
\begin{table}
\caption{GOE approximation. As for Table \protect\ref{tabgoe}, but using
$s_{(1)}^\prime$ defined using $\mu_N(\gamma)$ and $\tau_N(c)$
given by (\protect\ref{eqmodcentsc}), with $\gamma= 1/5$ and $c=1$}
\label{tabgoemod}
\begin{tabular*}{\textwidth}{@{\extracolsep{\fill}}lccccccccc@{}}
\hline
$\bolds{N+1}$&\textbf{0.01}&\textbf{0.05}&\textbf{0.1}&\textbf
{0.3}&\textbf{0.5}&\textbf{0.7}&\textbf{0.9}&\textbf{0.95}&\textbf
{0.99}\\
\hline
\phantom{00}2&0.022&0.073&0.127&0.319&0.505&0.696&0.897&0.950&0.991\\
\phantom{00}3&0.018&0.067&0.120&0.315&0.505&0.699&0.899&0.951&0.991\\
\phantom{00}4&0.017&0.063&0.116&0.311&0.505&0.699&0.900&0.951&0.991\\
\phantom{00}5&0.015&0.061&0.114&0.310&0.504&0.700&0.901&0.951&0.991\\
\phantom{0}10&0.013&0.056&0.107&0.305&0.502&0.700&0.901&0.951&0.991\\
\phantom{0}25&0.011&0.053&0.104&0.302&0.501&0.700&0.901&0.951&0.990\\
\phantom{0}50&0.011&0.052&0.102&0.301&0.500&0.699&0.900&0.950&0.990\\
\phantom{0}75&0.011&0.052&0.102&0.302&0.500&0.700&0.900&0.950&0.990\\
100&0.010&0.051&0.101&0.301&0.500&0.700&0.900&0.951&0.990\\
200&0.011&0.051&0.101&0.301&0.500&0.700&0.900&0.950&0.990\\
500&0.010&0.050&0.100&0.300&0.500&0.700&0.901&0.951&0.990\\
\hline
\end{tabular*}
%
\end{table}

%








\textit{Outline of proof}.
We use the operator norm convergence framework developed
in \citet{trwi04}; our focus, of course, is on achieving the second order
convergence rate results.
We use the Fredholm determinant representations for the finite
and limiting distribution functions in terms of the two-point
correlation kernels. A bound of Seiler--Simon, along with its
orthogonal case analog, bounds the difference in Fredholm determinants
in terms of the kernels. In turn, the kernels have integral
representations in terms of weighted Hermite polynomials, and so we
transfer bounds on convergence of Hermite polynomials to the Airy
function to bounds on the kernels and hence to bounds on the
probabilities.

Convenient uniform bounds on the convergence of weighted Hermite
polynomials to
the Airy function come from Liouville--Green theory,
which analyzes convergence of the solutions of the second-order
differential equation satisfied by the Hermite polynomials to those of
the equation for the Airy function.

The correlation kernels for finite $N$ involve polynomials of both
degrees~$N$ and $N-1$, each with its own Liouville--Green centering
$u_N$ and $u_{N-1}$.
The overall centering $\mu_N$ for the kernel and distribution function
must combine~$u_N$ and $u_{N-1}$ appropriately to ensure that the
generic $O(N^{-1/3})$ error terms cancel to uncover $O(N^{-2/3})$
convergence.
In the unitary case, simple averaging suffices: $\mu_N = (u_N +
u_{N-1})/2$. For the orthogonal setting, we use a~formula expressing
the $\operatorname{GOE}_{N+1}$ kernel in terms of the $\operatorname
{GUE}_N$ kernel plus a rank one
kernel, and obtain cancellation of $O(N^{-1/3})$ errors from these
two components.

The Hermite polynomial approximation results are summarized in
Section~\ref{secherm-polyn-asympt}. The unitary proof, in Section~\ref{secunitary-case}, is a necessary preparation
for the orthogonal case in Section~\ref{secorthogonal-case}.

\textit{Reproducible code:}
\textsc{Matlab} files to produce the figures and tables are available
at the second author's website.

\textit{Related work}. Convergence rate results at $O(N^{-2/3})$
were obtained by \citet{elka04p} for LUE, \citet{john08p} for JUE and
JOE, and \citet{ma11} for LOE. The study of Edgeworth-type expansions
for GUE
and LUE was initiated by Choup (\citeyear{chou06,chou08}), who noted that
$N^{-1/3}$ terms in these expansions can be removed by specific
choices of the centering constant.

The Tracy--Widom limit laws for the largest eigenvalue hold much more
generally---such universality results are an active subject of
research. For Hermitian Wigner matrices, see \citet{tavu10} and
references therein, and for covariance matrices \citet{sosh01a} and
\citet{pech07}.



\section{Hermite polynomial asymptotics}
\label{secherm-polyn-asympt}

The Hermite polynomials, $H_k(x)$ in notation of
\citeauthor{szeg67} [(\citeyear{szeg67}), Chapter 4], are orthogonal with respect to the weight
function $w(x) = e^{-x^2}$ on $(-\infty,\infty)$.
The ``oscillator wave functions'' are normalized, weighted versions
%
\[
\phi_k(x) = h_k^{-1/2} e^{-x^2/2} H_k(x),
\]
with $h_k = \int H_k^2(x) e^{-x^2} \,dx = \sqrt{\pi} 2^k k!$

Classical Plancherel--Rotach asymptotics for $H_N(x)$ near the largest
zero, \citeauthor{szeg67} [(\citeyear{szeg67}), page 201] and \citeauthor{agz10} [(\citeyear{agz10}), Section~3.7.2],
establish that, for $m_N = \sqrt{2N}$ and
$\tau_N = 2^{-1/2} N^{-1/6}$,
%
\begin{equation}
\label{eqplr}
(2N)^{1/4} \tau_N \phi_N (m_N + s \tau_N) \rightarrow A(s),
\end{equation}
where throughout we use $A$ to denote the Airy function $\mbox{\rm Ai}$.

We will need to explicitly bound the error in the convergence in
\eqref{eqplr}.
There is now a substantial literature on asymptotic approximations to
Hermite polynomials, using, for example,
the steepest descent method for integrals [e.g., \citet{shi08}],
the nonlinear steepest descent method for Riemann--Hilbert problems
[e.g., \citet{wozh07}]
and recurrence relations [e.g., \citet{wawo11}].
Much of this recent attention has focused on expansions for
$H_N(\sqrt{2N+1} \xi)$ and $\phi_N(\sqrt{2N+1} \xi)$ that are
valid uniformly for large regions of $\xi$.

For this work, however, we need
more detailed information for $\xi= 1 + \sigma_N s$ near $1$,
and specifically uniform bounds for the
error of Airy approximation for both $\phi_N$ \textit{and its
derivative} that have exponential decay in the variable $s$
and rate $N^{-2/3}$; cf. Proposition~\ref{prop-localbds} below.
We have not found this extra detail explicitly in the literature, and
since the Liouville--Green discussion of \citeauthor{olve74} [(\citeyear{olve74}), Chapter~11]
comes with
ready-made bounds for approximation error for both $\phi_N$ and
$\phi_N^\prime$, we use this as a starting point for extracting, in
the \hyperref[app]{Appendix}, the
specific bounds we need.
In this section, we explain just enough of the approach to describe
the bounds we need.

The Liouville--Green (LG) approach
relies on the fact that Hermite
polynomials, and hence $\phi_N$, satisfy a second order
differential equation,
%
\begin{equation}
\label{eqwNdiffeq}
\phi_N''(x) = \{ x^2 - (2N+1) \} \phi_N(x).
\end{equation}
Rescaling the $x$ axis via $x = \sqrt{2N+1} \xi$,
and setting $w_N(\xi) = \phi_N(x)$,
the equation takes the form
%
\begin{equation}
\label{eqhermite-de}
w_N''(\xi) =
\kappa_N^2 f(\xi) w_N(\xi),
\end{equation}
with
%
\[
\kappa_N = 2N+1,\qquad f(\xi) = \xi^2 - 1.
\]

The \textit{turning points} of the differential equation are the zeros
of $f$, namely $\xi_{\pm} = \pm1,$ so named because each separates an
interval in which the solution is of exponential type from one in
which the solution oscillates. The LG transformation
introduces new independent and dependent variables $\zeta$ and~$W$ via
the equations
%
\begin{equation}
\label{eqnewvarsa}
\zeta\biggl( \frac{d\zeta}{d\xi} \biggr)^2 = f(\xi),\qquad
W = \biggl( \frac{d\zeta}{d\xi} \biggr)^{1/2} w_N.
\end{equation}
More precisely, we take
%
\begin{equation}
\label{eqnewvarsb}
(2/3) \zeta^{3/2} (\xi) = \int_1^\xi\sqrt{f(\xi')} \,d\xi'.
\end{equation}
The transform $W$ approximately satisfies the Airy equation
$W''(\zeta) =\break \kappa_N^2 \zeta W(\zeta)$, which has
linearly independent solutions in terms of Airy functions,
traditionally denoted by $\operatorname{Ai}(\kappa^{2/3} \zeta)$ and
$\operatorname{Bi}(\kappa^{2/3} \zeta)$. Our interest lies in
approximating the
\textit{recessive} solution $\operatorname{Ai}(\kappa^{2/3} \zeta)$.

As described in more detail in the \hyperref[app]{Appendix},
the error in
the Liouville--Green approximation can be bounded, and one arrives at
%
\begin{equation}
\label{eqNm1error}
\bar\phi_N (x) = (2N)^{1/4} \tau_N \phi_N(x)
= \bar e_N r(\xi) \{ \Ai( \kappa_N^{2/3} \zeta) +
O(N^{-1}) \}.
\end{equation}
Here $r(\xi) = [ \dot\zeta(\xi)/ \dot\zeta(1) ]^{-1/2}$ is
approximately $1$ for $\xi$ near $1$, and $\bar e_N = 1 + O(N^{-1})$.
This is, then, a version of~\eqref{eqplr} with an error term of order
$O(N^{-1})$, but with Airy function argument $\kappa_N^{2/3} \zeta$
rather than $s$.

We focus on $x$ near $u_N = \sqrt{2N + 1}$, that is, on $\xi$ near the
upper turning point $\xi_+ = 1$.
Introduce the rescaled variable $s$ through $\xi= 1 + \sigma_N s$.
To more closely match the result~\eqref{eqplr}, we want $\sigma_N$ to
be chosen
so that the Airy function argument
%
\begin{equation}
\label{eqlingrowth}
\kappa_N^{2/3} \zeta(1+ \sigma_N s) \approx s
\end{equation}
for $s$ in a suitably large range.
A Taylor expansion of the left-hand side yields
%
\begin{equation}
\kappa_N^{2/3} \bigl(\zeta(1) + \sigma_N \dot\zeta(1) s
+ \tfrac{1}{2}\sigma_N^2 s^2 \ddot\zeta^*\bigr).
\end{equation}
Since $\zeta(1) = 0$ and $\dot\zeta(1) = 2^{1/3}$, as
follows from~\eqref{eqnewvarsa} and~\eqref{eqnewvarsb}, we obtain
\eqref{eqlingrowth} by any choice of the form
%
\[
\sigma_N = \tfrac{1}{2}N^{-2/3}\bigl (1 + o(1)\bigr).
\]
For such a choice,~\eqref{eqtaylor} shows that
%
\[
\kappa_N^{2/3} \zeta(1 + \sigma_N s) = s + O( s^2 N^{-2/3}).
\]
Thus, to replace $\Ai(\kappa_N^{2/3}
\zeta)$ in~\eqref{eqNm1error} by $\Ai(s)$ entails, in general,
accepting an error term of $O(N^{-2/3})$ instead of $O(N^{-1})$, and
so we use this error scale henceforth.

\begin{figure}

\includegraphics{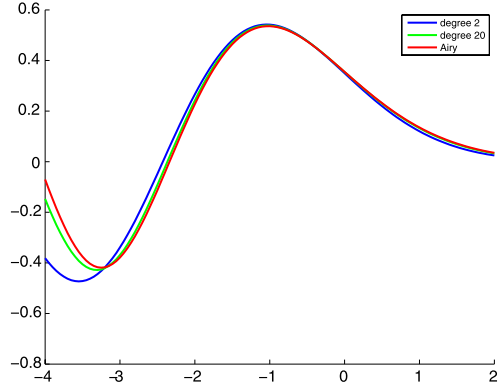}

\caption{Plots of $\bar\phi_N(u_N + \tau_N s)$ for $N = 2, 20$
compared with the Airy function $A(s)$.}
\label{figpolys}
\end{figure}

With the specific choice $\sigma_N = \tau_N/u_N$, we will show that
for $s \geq s_L$,
%
\[
|\bar\phi_N( u_N + s \tau_N) - \Ai(s) | \leq C N^{-2/3} e^{-s/2}.
\]
Figure~\ref{figpolys} shows that, for values of $s$ corresponding to
the bulk of the support of $F_1$, the approximation is tolerably good
even for $N=2$.

In fact,
since the two-point correlation functions depend on both $\phi_N$ and
$\phi_{N-1}$, we need such approximations both for
$\bar\phi_N$ and for $\bar\phi_{N-1}$.
For $\bar\phi_{N-1} = (2N)^{1/4} \tau_N \phi_{N-1}$, the corresponding turning point is at
$u_{N-1} = \sqrt{2N - 1}$, though we still use the same scale factor
$\tau_N$.
We use the notation $\bar{\phi}_{Nj}$, with $Nj = N$ or $N-1$,
respectively, to refer to both cases.
In addition, for the GOE case, bounds on the convergence of the
derivative is also required.
In the \hyperref[app]{Appendix}, we establish

\begin{proposition}
\label{prop-localbds}
Let $s_L \in\mathbb{R}$. 
For $s \geq s_L$, we have
%
\begin{eqnarray}
|\bar{ \phi}_{Nj}( u_{Nj} + s \tau_N)|
& \leq &C e^{-s}, \label{eqphiNjbd} \\
| \bar{ \phi}_{Nj} (u_{Nj} + s \tau_N) - \Ai(s)|
& \leq& C N^{-2/3} e^{-s/2}, \label{eqdelphiNj}
\end{eqnarray}
%
where the error bounds are uniform in $s \geq s_L$ and $N \geq N_0(s_L)$.
The same bounds hold, with modified constants $C$, when
$\bar\phi_{Nj}$ and $\Ai$ are replaced by
$\tau_N \bar\phi_{Nj}^\prime$ and $\Ai^{ \prime}$.
\end{proposition}

We record here also some corresponding exponential decay bounds for the Airy
function and its derivatives $\Ai^{(i)}$.
Indeed, given $s_L$, there exist
constants $C_i(s_L)$ such that
%
\begin{equation}
\label{eqAiFacts}
\bigl| s^k \Ai^{(i)}(s)\bigr| \leq C_i(s_L) e^{-s},\qquad s \geq s_L, i = 0,1,2.
\end{equation}

Proposition~\ref{prop-localbds} provides good Airy approximations for
both $\bar\phi_N$ and $\bar\phi_{N-1}$, but with differing centering
values, $u_N$ and $u_{N-1}$, respectively. To obtain scaling limits
for the correlation kernels in GUE and GOE, we need to combine these
centerings in a manner appropriate to each case.

It is convenient to express these centering shifts in the
rescaled variable~$s$.
Thus, set
%
\begin{eqnarray}
\label{eqkextension}
\phi_\tau(s; k) & =& \bar\phi_N\bigl(u_N + \tau_N(s + k \Delta_N)\bigr),
\nonumber
\\[-8pt]
\\[-8pt]
\nonumber
\psi_\tau(s; l) & =& \bar\phi_{N-1}\bigl(u_{N-1} + \tau_N(s + l \Delta_N)\bigr),
\end{eqnarray}
where
%
\[
\Delta_N = (u_N - u_{N-1})/\tau_N = N^{-1/3}\bigl( 1 + 2^{-5} N^{-2} +
O(N^{-4})\bigr)
\]
---indeed $2^{1/3} \Delta_2 = 1.0080!$
We obtain extensions of Proposition~\ref{prop-localbds}: indeed from~\eqref{eqphiNjbd},
%
\[
|\phi_\tau(s;k)| \leq C e^{-(s + k \Delta_N)} \leq C e^{-s},
\]
and a similar bound holds for $|\psi_\tau(s;k)|$.
From (\ref{eqkextension}) and Proposition~\ref{prop-localbds},
%
\begin{eqnarray}\label{eqtwo-term}
\phi_\tau(s;k) & =& \Ai(s + k \Delta_N) + O(N^{-2/3} e^{-s/2})
\nonumber
\\[-8pt]
\\[-8pt]
\nonumber
& = &\Ai(s) + k \Delta_N \Ai'(s) + O(N^{-2/3} e^{-s/2}),
\end{eqnarray}
since $\frac{1}{2}(k \Delta_N)^2 | \Ai''(s^*)| \leq C N^{-2/3}
e^{-s}$ using
(\ref{eqAiFacts}) and $|s^* - s| \leq k \Delta_N$.

More generally, but by identical arguments, for $r = 0,1$ we have
%
\begin{equation}
\label{eqkeffect}
\phi_\tau^{(r)}(s;k)
= \Ai^{(r)}(s) + k \Delta_N \Ai^{(r+1)}(s) + O(N^{-2/3} e^{-s/2}),
\end{equation}
and correspondingly for $\psi_\tau^{(r)}(s;k)$.
As a byproduct of a steepest descent analysis for Laguerre
polynomials, Choup (\citeyear{chou06,chou08})
derived a three-term asymptotic expansion for $\phi_N$ and
$\phi_{N-1}$ whose first two terms agree with~(\ref{eqtwo-term}),
though without the uniform error bounds in $N$ and $s$ of
Proposition~\ref{prop-localbds}.


\section{Unitary case}
\label{secunitary-case}
\mbox{}
\begin{pf*}{Proof of Theorem \protect\ref{thunitary-main-res}}
The argument has three components: first we recall determinantal
representations of the eigenvalue probabilities $F_{N,2}(x_0)$ and
limiting value $F_2(s_0)$, along with integral representations of the
associated correlation kernels.
Second we set up the rescaling that connects $x_0$ and $s_0$, and
finally establish the convergence bounds.

The two-point correlation kernel for $\operatorname{GUE}_N$
%
\[
S_{N,2} (x,y) = \sum_{k=0}^{N-1} \phi_k(x) \phi_k(y)
\]
has a useful integral representation [\citet{trwi96}, equation~(57)].
Set\setcounter{footnote}{1}\footnote{Note: our definitions differ by a factor
$\sqrt2$ from those of Tracy and Widom.}
%
\begin{equation}
\label{eqphi-N-hat-S-N-hat0}
\phi(x) = (2N)^{1/4} \phi_N(x),\qquad
\psi(x) = (2N)^{1/4} \phi_{N-1}(x),
\end{equation}
then
%
\begin{equation}
\label{eqSN2-integ}
S_{N,2} (x,y)
= \frac{1}{2} \int_0^\infty[\phi(x+z) \psi(y+z) + \psi(x+z)
\phi(y+z)] \,dz.
\end{equation}

The distribution of $x_{(1)}$ may be expressed as a Fredholm
determinant
%
\begin{equation}
F_{N,2}(x_0) = P \Bigl\{ \max_{1 \leq k \leq N} x_k \leq x_0 \Bigr\} =
\det(I - S_{N,2} \chi_0 ),
\label{eqdetlargest}
\end{equation}
where $\chi_0(x) = I_{(x_0, \infty)}(x)$ and the operator $S_{N,2}
\chi_0$ is defined via
%
\[
(S_{N,2} \chi_0) g(x) = \int_{x_0}^\infty S_{N,2}(x,y) g(y) \,dy.
\]
Equivalently, we may speak of $S_{N,2}$ as an operator on
$L_2(x_0,\infty)$ with kernel $S_{N,2}(x,y)$.
On this understanding, we drop further explicit reference to $\chi_0$.

Now change variables, setting $x = \tau(s) = \mu_N +
\tau_N s,$ with $\mu_N$ yet to be determined, and
$x_0 = \tau(s_0)$.
Set also
%
\begin{equation}
\label{eqStaudef}
S_\tau(s,t) = \tau_N S_{N,2}(\mu_N + \tau_N s, \mu_N + \tau_N t).\vadjust{\goodbreak}
\end{equation}
Defining
%
\begin{equation}
\label{eqphipsitau0}
\phi_\tau(s) = \tau_N \phi(\mu_N + \tau_N s),\qquad
\psi_\tau(s) = \tau_N \psi(\mu_N + \tau_N s),
\end{equation}
it is clear that (\ref{eqSN2-integ}) becomes
%
\begin{equation}
S_\tau(s,t) = \frac{1}{2} \int_0^\infty
[\phi_\tau(s+z) \psi_\tau(t+z) + \psi_\tau(s+z) \phi_\tau(t+z)]
\,dz. \label{eqStau}
\end{equation}

Since $S_{N,2}$ and $S_\tau$ have
the same eigenvalues, $\det(I-S_{N,2}) = \det(I-S_\tau)$, and so
%
\begin{equation}
\label{eqN-prob}
P \{ (\max x_k - \mu_N)/\tau_N \leq s_0 \} = \det(I-S_\tau).
\end{equation}

\citet{trwi94} showed that the limiting distribution $F_2$
also has a determinantal representation
%
\[
F_2(s_0) = \det(I - S_A),
\]
where $S_A$ denotes the Airy operator on $L^2(s_0,\infty)$ with the kernel
having the form
%
\begin{equation}
\label{eqSbar}
S_A (s,t) = \int_0^\infty\Ai(s+z) \Ai(t+z) \,dz.
\end{equation}
To derive bounds on the convergence of $F_{N,2}(x_0)$ to $F_2(s_0)$,
we use a bound due to \citet{sesi75},
%
\begin{equation} \label{eqsei-sim}
| \det(I-S_{\tau}) - \det(I- S_A)|
\leq\| S_\tau- S_A \|_1 \exp( \| S_\tau\|_1 + \| S_A
\|_1 + 1 ).
\end{equation}
%
Here $\| \cdot\|_1$ denotes trace class norm on
operators on $L^2(s_0,\infty)$.
This bound reduces the convergence
question to study of convergence of the kernel $S_\tau(s,t)$ to
$S_A(s,t)$.

Given functions $a$ and $b$, denote by $a \diamond b$ the operator
having kernel
%
\[
(a \diamond b)(s,t) = \int_{0}^\infty a(s+z) b(t+z) \,dz.
\]
In this notation, the kernel difference becomes
%
\[
S_\tau- S_A
= \tfrac{1}{2}( \phi_\tau\diamond\psi_\tau+ \psi_\tau\diamond
\phi_\tau)
- A \diamond A.
\]
%
To facilitate convergence arguments, we rewrite this as
\begin{eqnarray*}
8(S_\tau- S_A)
& =& (\phi_\tau+ \psi_\tau+ 2 A) \diamond(\phi_\tau+ \psi_\tau-
2A) \\
&&{} + (\phi_\tau+ \psi_\tau- 2 A) \diamond(\phi_\tau+ \psi_\tau+
2A) - (\phi_\tau- \psi_\tau) \diamond(\phi_\tau- \psi_\tau).
\end{eqnarray*}
%

Recall that the centering constant $\mu_N$ was left unspecified in the
definitions of $\phi_\tau$ and $\psi_\tau$ in
(\ref{eqphipsitau0}). We now choose $\mu_N$ so that each term in the
preceding decomposition is $O(N^{-2/3})$. This amounts to choosing the
shifts $k$ and~$l$ in (\ref{eqkextension}) to satisfy two
constraints. First, the centerings $\mu_N = u_N + k \tau_N \Delta_N$
and $\mu_N = u_{N-1} + l \tau_N \Delta_N$ must agree,
so that necessarily $l = k+1$.
Second, the\vadjust{\goodbreak} $N^{-1/3}$ term must drop out in the expansion for $\phi
_\tau+
\psi_\tau$ given by (\ref{eqkeffect}), so that $l = -k$.
We therefore must have, for the present unitary case,
%
\begin{eqnarray}
\label{eqphipsireal}
\phi_\tau(s)
& =& \phi_\tau\bigl(s ; - \tfrac{1}{2}\bigr),
\nonumber
\\[-8pt]
\\[-8pt]
\nonumber
\psi_\tau(s)
& =& \psi_\tau\bigl(s ; \tfrac{1}{2}\bigr)
\end{eqnarray}
which entails that $\mu_N = (u_N + u_{N-1})/2$ as
was used in Table~\ref{tabguemod}.
From Proposition~\ref{prop-localbds} and the succeeding discussion we
obtain

\begin{corollary}[(Complex Case)]
\label{lemkeybounds-com}
Let $\phi_\tau$ and $\psi_\tau$ be defined by (\ref{eqphipsireal})
and~(\ref{eqkextension}).
Given $s_L \in\mathbb{R}$, there exists $C = C(s_L)$
such that for $N \geq N(s_L)$ and $s \geq s_L$,
%
\begin{eqnarray}
| \phi_\tau(s) | , | \psi_\tau(s) | & \leq& C
e^{-s}, \label{eqphi-global-bound}\\
| \phi_\tau(s) - A(s) |, | \psi_\tau(s) -
A (s) | & \leq& C N^{-1/3} e^{-s/2},
\label{eqphitau} \\
| \phi_\tau(s) + \psi_\tau(s) - 2 A (s) | & \leq& C N^{-2/3}
e^{-s/2}. \label{eqphipsitau}
\end{eqnarray}
\end{corollary}

We will need some simple bounds for certain norms of $a \diamond b$.
In the unitary case, we need the trace norm of $a \diamond b$ as an
operator on $L^2(s_0, \infty)$.
In the orthogonal case, we need
the weighted $L^2$-spaces
$L^2((s_0,\infty), \rho(s) \,ds)$ and
$L^2((s_0,\infty), \rho^{-1}(s) \,ds)$ for a
weight function $\rho$ such that the reciprocal $\rho^{-1} \in
L^1(\mathbb{R})$.
Further details are given in Section~\ref{secorthogonal-case}.
For some $\gamma\geq0$, let
%
\begin{equation}
\label{eqrhospec}
\rho(s) = e^{\gamma|s|}.
\end{equation}
In this section $\gamma= 0$, while values of $\gamma> 0$ will be specified
later for GOE.

\begin{proposition}
\label{propop-norm}
Let weight functions $\rho_1, \rho_2$ be chosen from $\{
\rho, 1/\rho\}$, where~$\rho$ is given by (\ref{eqrhospec}), and
consider the Hilbert--Schmidt norm of operator
$a \diamond b\dvtx  L^2(\rho_2) \rightarrow L^2(\rho_1)$.
Assume that, for $s \geq s_0$,
%
\[
|a(s)| \leq a_N e^{-a_1 s},\qquad
|b(s)| \leq b_N e^{-b_1 s}.
\]
If $0 \leq\gamma< 2 \min(a_1, b_1)$, then
%
\begin{equation}
\label{eqop-bound}
\| a \diamond b \|_{HS} \leq C \frac{a_N b_N}{a_1 + b_1} e^{-(a_1 +
b_1) s_0 \pm\gamma|s_0|},
\end{equation}
where $C = C(a_1, b_1, \gamma) = [(a_1 - \gamma/2)(b_1 -
\gamma/2)]^{-1/2}.$
If $\rho_1 = \rho_2,$ then the trace norm satisfies the same bound.
\end{proposition}

This is a special case of \citeauthor{john08p} [(\citeyear{john08p}), Lemma 7].
In the present unitary case,
we apply Proposition~\ref{propop-norm}, with $\gamma= 0$, to bound the
trace norm of each term on the right-hand side, using
(\ref{eqphi-global-bound})--(\ref{eqphipsitau}).
For each of the three terms, we find that
$a_1 + b_1 \geq1$ and $a_N b_N \leq
C N^{-2/3}$, so that
%
\[
\| S_\tau- S_A \|_1 \leq C N^{-2/3} e^{-s_0}.
\]
We may similarly conclude that
%
\[
\| S_\tau\|_1 \leq C e^{-2 s_0},\qquad
\| S_A \|_1 \leq C e^{-2 s_0}.\vadjust{\goodbreak}
\]
Indeed, bounds~\eqref{eqphi-global-bound} and~\eqref{eqAiFacts}
show that in each case, $a_N$ and $b_N \leq C(s_0)$, and that $a_1 =
b_1 = 1$.

Combining the two previous displays with the Seiler--Simon bound
\eqref{eqsei-sim}, we obtain Theorem~\ref{thunitary-main-res}.
\end{pf*}

\section{Orthogonal case}
\label{secorthogonal-case}

To establish Theorem~\ref{thorthogonal-main-res}, we again follow the
outline of proof given in Section~\ref{secintroduction}.

$1^\circ.$ Assume that $N+1$ is even.
\citet{trwi98} gave a
derivation\footnote{\citet{sinc09} extended Tracy and Widom's
derivation to cover $N + 1$ odd, but we do not pursue this here. See
also \citet{foma09}.}
of the determinant representation
%
\begin{equation}
\label{eqOrthogProb0}
P \Bigl\{ \max_{1 \leq k \leq N+1} x_k \leq x_0 \Bigr\} =
\sqrt{ \det(I - K_{N+1} \chi_0 ) }.
\end{equation}
Here $K_{N+1}$ is a $2 \times2$-matrix valued operator
%
\begin{equation}
\label{eqmatrix-rep0}
K_{N+1}(x,y) = (L S_{N+1,1})(x,y) + K^\varepsilon(x,y),
\end{equation}
where
%
\begin{equation}
\label{eqKeps-def}
L =
\pmatrix{
I & - \partial_2 \vspace*{2pt}\cr \varepsilon_1 & T},
\qquad
K^\varepsilon(x,y) =
\pmatrix{
0 & 0 \vspace*{2pt}\cr -\varepsilon(x-y) & 0}
.
\end{equation}
Here $\partial_2$ denotes the operator of partial differentiation with
respect to the second variable, and $\varepsilon_1$ the operator of
convolution in the first variable with the function $\varepsilon(x) =
\frac{1}{2}
\operatorname{sgn} (x).$
Thus $(\varepsilon_1 S)(x,y) = \int\varepsilon(x-u) S(u,y) \,du.$
Finally~$T$ denotes transposition of variables $TS(x,y) = S(y,x).$
The scalar kernel
%
\[
S_{N+1,1}(x,y)
= \sum_{n=0}^N \phi_n(x) \phi_n(y) + \sqrt{ \frac{N+1}{2}}
\phi_N(x) \varepsilon\phi_{N+1}(y),
\]
and \citet{afnm00} observe that it may be rewritten as
%
\begin{equation}
\label{eqRealKernel0}
S_{N+1,1}(x,y)
=S_{N,2}(x,y) + \tfrac{1}{2}\phi(x) \varepsilon\psi(y),
\end{equation}
where $\phi$ and $\psi$ are as defined at
(\ref{eqphi-N-hat-S-N-hat0}).
The orthogonal kernel is thus expressed in terms of the unitary kernel
and a rank one remainder term.
The formula allows convergence results from the unitary case to be
reused, with relatively minor modification.

$2^\circ.$
The limiting distribution has a corresponding determinantal
representation
%
\[
F_1(s_0) = \sqrt{ \det(I-K_{\mathrm{GOE}}) }.
\]
To state the \citet{trwi04} form for $K_{\mathrm{GOE}}$, and for the
convergence argument to follow, it is helpful to rewrite expressions
involving~$\varepsilon$ in terms\vadjust{\goodbreak} of the right-tail integration operator
$ (\tilde\varepsilon g)(s) = \int_s^\infty g(u) \,du$
and for kernels $A(s,t)$ in the form
$ (\tilde\varepsilon_1 A)(s,t) = \int_s^\infty A(u,t) \,du.$
This is due to the oscillatory behavior of the
Airy function in the \textit{left} tail.
We write $A \otimes B$ for the operator whose kernel is $A(s)
B(t)$.
The Tracy--Widom expression states that
%
\begin{equation}
\label{eqGOEker}
K_{\mathrm{GOE}}(s,t) =
\pmatrix{
S(s,t) & SD(s,t) \vspace*{2pt}\cr
IS(s,t) - \varepsilon(s-t) & S(t,s)},
\end{equation}
and the entries of $K_{\mathrm{GOE}}$ are given by
%
\begin{eqnarray} \label{eqKGOEform-1}
S(s,t) & =& \bigl(S_A - \tfrac{1}{2}A \otimes\tilde\varepsilon A\bigr)(s,t) +
\tfrac{1}{2}A(s),\nonumber \\
SD(s,t) & =& - \partial_2 \bigl(S_A - \tfrac{1}{2}A \otimes\tilde
\varepsilon
A\bigr)(s,t), \\
IS(s,t) & = &- \tilde\varepsilon_1 \bigl(S_A - \tfrac{1}{2}A \otimes\tilde
\varepsilon A\bigr)(s,t) - \tfrac{1}{2}(\tilde\varepsilon A)(s) + \tfrac
{1}{2}(\tilde
\varepsilon A)(t),\nonumber
\end{eqnarray}
where $S_A$ is the Airy kernel defined at~\eqref{eqSbar}.

Defining operator matrices
%
\[
\tilde L =
\pmatrix{ I & - \partial_2 \vspace*{2pt}\cr - \tilde\varepsilon_1 & T}
,\qquad
L_1 =
\pmatrix{
I & 0 \vspace*{2pt}\cr - \tilde\varepsilon& 0},\qquad
L_2 =
\pmatrix{
0 & 0 \vspace*{2pt}\cr \tilde\varepsilon& I},
\]
we may rewrite (\ref{eqGOEker}) in the form
%
\begin{equation}
\label{eqKGOE}
K_{\mathrm{GOE}} =
\tilde L\bigl(S_A - \tfrac{1}{2}A \otimes\tilde\varepsilon A\bigr)
+ \tfrac{1}{2}L_1 A(s) + \tfrac{1}{2}L_2 A(t) + K^\varepsilon.
\end{equation}

$3^\circ$. We turn to a linear rescaling of
formulas (\ref{eqOrthogProb0}) and
(\ref{eqmatrix-rep0}).
We again set $x = \tau(s) =
\mu_N + \tau_N s$ and $ y = \tau(t) = \mu_N + \tau_N t$, but now
with $\mu_N = \mu_N^R$ to be determined anew in this orthogonal case;
see $4^\circ$ below.
Define $\phi_\tau$ and $\psi_\tau$ as before by (\ref{eqphipsitau0});
we occasionally write $\phi_\tau^R$ and $\psi_\tau^R$ to emphasize the
different centering. We have
%
\begin{eqnarray}\label{eqSRStau}
S_\tau^R(s, t) & :=&
\tau_N S_{N+1,1}( \tau(s), \tau(t)) \nonumber\\
& = &\tau_N S_{N,2}( \tau(s), \tau(t)) + \tfrac{1}{2}\tau_N \phi
(\tau(s))
(\varepsilon\psi)(\tau(t)) \\
& =& S_\tau(s,t) + \tfrac{1}{2}\phi_\tau(s) (\varepsilon\psi_\tau
)(t),\nonumber
\end{eqnarray}
where we have used $(\varepsilon\psi) (y) = (\varepsilon\psi_\tau) (t)$
for a linear rescaling.

Now $\det(I - K_{N+1} \chi_0) = \det(I - \bar K_\tau)$, where
\begin{eqnarray*}
\bar K_\tau(s,t) & =& \tau_N K_{N+1}(\tau(s), \tau(t)) \\
& =& \tau_N
\pmatrix{
I & - \partial_2 \vspace*{2pt}\cr \varepsilon_1 & T}
S_{N+1,1}(\tau(s),\tau(t)) + \tau_N K^\varepsilon( \tau(s), \tau(t))
\\
& =&
\pmatrix{
I & -\tau_N^{-1} \partial_2 \vspace*{2pt}\cr
\tau_N \varepsilon_1 & T}
S_\tau^R(s,t) + \tau_N K^\varepsilon(s,t),
\end{eqnarray*}
where $K^\varepsilon$ was defined at (\ref{eqKeps-def}).
Since $\det(I- \bar K_\tau)$ is unchanged if the lower left entry is
divided by $\tau_N$ and the upper right entry multiplied by
$\tau_N$,
%
\begin{equation}
\label{eqscalerep-a}
\det(I - K_{N+1} \chi_0) = \det(I - K_\tau),
\end{equation}
where $K_\tau$ is an operator with matrix kernel
%
\begin{equation}
\label{eqscalerep-b}
K_\tau(s,t) = (L S_\tau^R)(s,t) + K^\varepsilon(s,t).
\end{equation}

Now we rewrite $LS_\tau^R$ using $\tilde\varepsilon$ and $\tilde
\varepsilon_1$. First, define
%
\begin{equation}
\label{eqbetaN-1def}
\beta_{N-1}= \frac{1}{2} \int_{-\infty}^\infty\psi_\tau
= \frac{1}{2} (2N)^{1/4} \int_{-\infty}^\infty\phi_{N-1},
\end{equation}
and observe that
$\varepsilon\psi_\tau= \beta_{N+1} - \tilde\varepsilon\psi_\tau$.
Thus
%
\[
L S_\tau^R = L\bigl(S_\tau- \tfrac{1}{2}\phi_\tau\otimes\tilde
\varepsilon
\psi_\tau\bigr) + \tfrac{1}{2}\beta_{N-1} L( \phi_\tau\otimes1).
\]
Now $L = \tilde L +
\left({{0 \atop \varepsilon_1 + \tilde\varepsilon_1} \enskip
{0 \atop 0}}\right)$
and $2(\varepsilon_1 + \tilde\varepsilon_1)$ amounts to integration over
$\mathbb{R}$ in the first slot.
From (\ref{eqStau}), after interchanging orders of integration and
using $\int\phi_\tau= 0$, we obtain
%
\[
\int_{-\infty}^\infty S_\tau(s,t) \,ds
= \int_0^\infty\beta_{N-1} \phi_\tau(t+z) \,dz
= \beta_{N-1} \tilde\varepsilon\phi_\tau(t),
\]
and then
%
\[
[ (L - \tilde L)(S_\tau- \phi_\tau\otimes\tilde\varepsilon
\psi_\tau) ]_{2,1} = \tfrac{1}{2}\beta_{N-1} \otimes\tilde
\varepsilon\phi_\tau
\]
as the only nonzero entry of the matrix on the left-hand side.
Combining the last two displays with~\eqref{eqscalerep-b}, we get
%
\begin{equation}
\label{eqLSRT}
K_\tau= \tilde L\bigl(S_\tau- \tfrac{1}{2}\phi_\tau\otimes\tilde
\varepsilon
\psi_\tau\bigr) + \tfrac{1}{2}\beta_{N-1} [L_1 \phi_\tau(s) + L_2 \phi
_\tau(t)]
+ K^\varepsilon.
\end{equation}


$4^\circ.$
We now look at the $(1,1)$ terms in~\eqref{eqKGOE} and
\eqref{eqLSRT} in order to see,
somewhat informally, how the choice $\mu_N^R = u_N$ leads to
$O(N^{-2/3})$ convergence.
Thus, we examine the difference
%
\begin{equation}
\label{eq11diff}
\bigl[S_\tau- \tfrac{1}{2}\phi_\tau\otimes\tilde\varepsilon_\tau\bigr]
- \bigl[S_A - \tfrac{1}{2}A \otimes\tilde\varepsilon A\bigr]
+ \tfrac{1}{2}[ \beta_{N-1} \phi_\tau- A].
\end{equation}
From definitions
(\ref{eqkextension}) and expansions (\ref{eqkeffect}), this choice of
$\mu_N^R$ corresponds to
%
\begin{eqnarray}
\label{eqrealcase}
\phi_\tau^R (s) & =& \phi_\tau(s; 0)
= \A(s) + O(N^{-2/3}),
\nonumber
\\[-8pt]
\\[-8pt]
\nonumber
\psi_\tau^R (s) & =& \psi_\tau(s; 1)
= \A(s) + \Delta_N \A'(s) + O(N^{-2/3}).
\end{eqnarray}

We write $A_N = A + \Delta_N A^\prime$ and define
%
\[
S_{A_N} = \tfrac{1}{2}( A \diamond A_N + A_N \diamond A).
\]
From representation~\eqref{eqStau} and~\eqref{eqrealcase},
$ S_\tau= S_{A_N} + O(N^{-2/3}),$
while the identity
%
\[
S_{A_N} = S_A - \tfrac{1}{2}\Delta_N A \otimes A
\]
follows from
%
\[
( A \diamond A' + A' \diamond A)(s,t)
= \int_0^\infty\frac{d}{dz} [A(s+z) A(t+z) ] \,dz
= - A(s) A(t).
\]
Thus
%
\[
S_\tau= S_A - \tfrac{1}{2}\Delta_N A \otimes A + O(N^{-2/3}).
\]

Since $\tilde\varepsilon A^\prime= - A$, we have
$\tilde\varepsilon A_N = \tilde\varepsilon A - \Delta_N A$, and so
%
\[
\tfrac{1}{2}\phi_\tau\otimes\tilde\varepsilon\psi_\tau
= \tfrac{1}{2}A \otimes\tilde\varepsilon A - \tfrac{1}{2}\Delta_N A
\otimes A + O(N^{-2/3}).
\]
Forming the difference of the last two displays, we see an important
cancellation of the $O(N^{-1/3})$ terms involving $\Delta_N$, and
hence that the first two terms of~\eqref{eq11diff} together are
$O(N^{-2/3})$.

A computation with the recursion relation for Hermite polynomials
and then Stirling's formula [with its $O(N^{-1})$ error term] shows
that, as $N \rightarrow\infty,$
%
\[
\beta_{N-1} = \biggl( \frac{ \pi N}{2} \biggr)^{1/4}
\frac{ \sqrt{(N-1)!}}{2^{(N-1)/2} ((N-1)/2)!}
= 1 + O(N^{-1}).
\]
From this and~\eqref{eqrealcase}, it follows that the final term of
\eqref{eq11diff} is also $O(N^{-2/3})$.

$5^\circ$. To prepare for the convergence argument for the $2 \times
2$ matrix kernels, we combine~\eqref{eqKGOE} and~\eqref{eqLSRT}.
Noting also from our considerations above that
%
\[
S_A - \tfrac{1}{2}A \otimes\tilde\varepsilon A
= S_{A_N} - \tfrac{1}{2}A \otimes\tilde\varepsilon A_N,
\]
we obtain
the basic difference representation
%
\begin{eqnarray}
\label{eqdiff-rep}
K_\tau- K_{\mathrm{GOE}} & =&
\tilde L(S_\tau- S_{A_N}) - \tfrac{1}{2}\tilde L( \phi_\tau\otimes
\tilde
\varepsilon\psi_\tau- A \otimes\tilde\varepsilon A_N)
\nonumber
\\[-8pt]
\\[-8pt]
\nonumber
&&{} + \tfrac{1}{2}L_1[ \beta_{N-1} \phi_\tau- A](s)
+ \tfrac{1}{2}L_2[ \beta_{N-1} \phi_\tau- A](t),
\end{eqnarray}
from which we may expect to show $O(N^{-2/3})$ convergence, in
view of the fact that $\phi_\tau, \psi_\tau, S_\tau$ and $\beta_{N-1}$
merge, respectively, with $A, A_N, S_{A_N}$ and $1$ at rates of at least
$O(N^{-2/3})$.

$ 6^\circ.$ We now turn to study the convergence of
%
\begin{equation}
\label{eqFNp1s0}
F_{N+1,1}(s_0) = P \bigl\{ \bigl(x_{(1)} - \mu_N\bigr)/\tau_N \leq s_0 \bigr\}
= \sqrt{ \det(I-K_\tau) }
\end{equation}
to $F_1(s_0) = \sqrt{ \det(I - K_{\mathrm{GOE}})}$.
\citet{trwi04} describe with some care the nature of the operator
convergence of $K_{N+1} $ to $K_{\mathrm{GOE}}$ for the Gaussian finite
$N$ ensemble.
We adopt their framework of weighted $L^2$ spaces and regularized
$2$-determinants.
Thus, let $\rho$ be a weight function such that $\rho^{-1} \in
L^1(\mathbb{R})$ and all $\phi_N \in L^2(\rho)$.
Write $L^2(\rho)$ and $L^2(\rho^{-1})$ for the spaces
$L^2((s_0,\infty), \rho(s) \,ds)$ and
$L^2((s_0,\infty), \rho^{-1}(s) \,ds)$, respectively.

We consider $K_\tau$ and $K_{\mathrm{GOE}}$ as members of the collection
$\mathcal{B}$ of $2 \times2$ Hilbert--Schmidt operator matrices $B =
(B_{ij}, i,j = 1,2 )$ on $L^2(\rho) \oplus L^2(\rho^{-1})$ whose
diagonal entries are trace class. Note that $\varepsilon\dvtx  L^2(\rho)
\rightarrow L^2(\rho^{-1})$ as a consequence of the assumption that
$\rho^{-1} \in L^1.$
The specific $\rho$ that we use is defined in~\eqref{eqrhospec} with
$\gamma> 0.$

To analyze the convergence of $p_{N+1} = F_{N+1,1}(s_0)$ to $p_\infty=
F_1(s_0)$, we note that their difference is bounded by
$|p_{N+1}^2 - p_\infty^2|/ p_\infty$,
so that we are led to the difference of determinants
%
\begin{equation}
\label{eqdetdiff}
|F_{N+1,1}(s_0) - F(s_0)| \leq C(s_0) | \det(I-K_\tau) - \det(I-K_{\mathrm{GOE}})|.
\end{equation}

A Seiler--Simon-type bound on the matrix operator
determinant for operators in~$\mathcal{B}$ is established in
\citet{john08p}.\vspace*{-2pt}

\begin{proposition}
\label{propdet-diff}
For $B, B' \in\mathcal{B}$, we have
%
\[
| \det(I-B) - \det(I-B')|
\leq C(B,B') \Delta(B - B'),
\]
where
%
\[
\Delta(B) = \sum_{i=1}^2 \| B_{ii} \|_1 +
\sum_{i \neq j} \| B_{ij} \|_2.
\]
The coefficient has the form $C(B,B') = \sum_{j=1}^2 c_{1j}(
\operatorname{tr} B,
\operatorname{tr} B') c_{2j}(B,B')$, where $c_{1j}$ and $c_{2j}$ are
continuous functions, the latter with respect to the strong
(Hilbert--Schmidt norm) topology.\vspace*{-2pt}
\end{proposition}

Insert the conclusion of Proposition~\ref{propdet-diff} into
(\ref{eqdetdiff}) to obtain
%
\begin{equation}
\label{eqF1summary}
|F_{N+1,1}(s_0) - F_1(s_0)|
\leq C(s_0, K_\tau, K_{\mathrm{GOE}}) \Delta( K_\tau- K_{\mathrm{GOE}}).
\end{equation}
%
We exploit decomposition~\eqref{eqdiff-rep}, which we write in the form
%
\[
K_\tau- K_{\mathrm{GOE}} = \delta^I + \delta_0^F + \delta_1^F + \delta_2^F
\]
to distinguish a term involving integral kernels, $\delta^I = \tilde
L(S_\tau- S_{A_N})$ from terms involving finite rank operators.
We establish trace norm bounds for the diagonal
elements and Hilbert--Schmidt bounds for the off-diagonal entries. The
distinction between the two norms
is moot for the finite rank terms $\delta_i^F,$
so the trace bounds
are actually also needed only for the $\delta^I$ term.

For each term, we show $\| \delta_{ij} \| \leq C N^{-2/3}$, so that
$\Delta(K_\tau- K_{\mathrm{GOE}})$ is bounded above by $C N^{-2/3}$.
We have both $\| K_\tau- K_{\mathrm{GOE}} \|_2$ and $\operatorname{tr} K_\tau-
\operatorname{tr} K_{\mathrm{GOE}} $ converging to $0$ at $O(N^{-2/3})$ rate,
so that
$C(K_\tau, K_{\mathrm{GOE}})$ remains bounded as $N \to\infty.$

$7^\circ.$ To bound each term $\delta_{ij}$,
we need orthogonal case analogs of the uniform bounds
of Corollary~\ref{lemkeybounds-com}, but now for
$\phi_\tau^R, \psi_\tau^R$ and their integrals and derivatives.
From Proposition
\ref{prop-localbds} and the succeeding discussion, we obtain\vspace*{-2pt}

\begin{corollary}[(Real case)]
\label{cororthogonal-case-1}
Let $\phi_\tau$ and $\psi_\tau$ be defined by (\ref{eqrealcase})
and (\ref{eqkextension}).
Given $s_L \in\mathbb{R}$, there exists $C = C(s_L)$
such that for $N \geq N(s_L)$ and $s \geq s_L$,
%
\begin{eqnarray}
| \phi_\tau(s) |
& \leq& C e^{-s}, \label{eqphitau-realbd} \\
| \psi_\tau(s) |
& \leq& C e^{-s}, \label{eqpsitau-realbd} \\
| \phi_\tau(s) - \Ai(s) |
& \leq& C N^{-2/3} e^{-s/2},
\label{eqphitau-real} \\
| \psi_\tau(s) - \Ai(s) - \Delta_N \Ai^{ \prime} (s) |
& \leq& C N^{-2/3} e^{-s/2}.
\label{eqpsitau-real}
\end{eqnarray}
The same bounds hold, with modified constants $C$, when $\phi_\tau,
\psi_\tau, \Ai$ and $\Ai^{ \prime}$ are replaced, respectively, by
$\phi_\tau^\prime, \psi_\tau^\prime, \Ai^{ \prime}$ and
$\Ai^{ \prime\prime}$,
or when $\psi_\tau, A$ and $A^\prime$ are replaced by
$\tilde\varepsilon\psi_\tau, \tilde\varepsilon A$ and $\tilde
\varepsilon
A^\prime$.
\end{corollary}

$\delta^I$ \textit{term}. For $\delta^{I} =\tilde L[ S_\tau- S_{A_N}]$,
we use Proposition~\ref{propop-norm} 
to establish the needed Hilbert--Schmidt and trace norm bounds
for each entry in the $2 \times2$ matrix. We write
%
\[
S_\tau- S_{A_N} =
(\phi_\tau- A) \diamond\psi_\tau+
A \diamond(\psi_\tau- A_N) +
(\psi_\tau- A_N) \diamond\phi_\tau+
A_N \diamond(\phi_\tau- A).
\]
In turn, for $\partial_2( \bar S_\tau- S_{A_N})$ we replace the
second slot arguments $\psi_\tau, (\psi_\tau- A_N)$,
etc., by their derivatives, and for
$\tilde\varepsilon(\bar S_\tau- S_{A_N})$, we replace the first slot
arguments $(\phi_\tau- A)$,
etc., by their right tail integrals.

Consider, for example, the first term $(\phi_\tau- A) \diamond
\psi_\tau$.
We apply Proposition~\ref{propop-norm}
using~\eqref{eqphitau-realbd} and~\eqref{eqphitau-real} to set
%
\[
a_N = CN^{-2/3},\qquad b_N = C,\qquad a_1 = \tfrac{1}{2}, b_1 = 1.
\]
The argument is entirely parallel when $\partial_2$ and $\tilde
\varepsilon_1$ is applied to $(\phi_\tau- A) \diamond
\psi_\tau$, and also for each of the second through
fourth terms.
Thus, if $D_{ij}$ denotes any matrix entry in any component of
$\delta^{I}$, we obtain
%
\begin{equation} \label{eqdelRI}
\| D_{ij} \| \leq C N^{-2/3} e^{-3s_0/2 + \gamma|s_0|}.
\end{equation}

\textit{Finite rank terms}.
As \citet{trwi04} note, the norm of a rank-one kernel $u(x) v(y)$, when
regarded as an operator $u \otimes v$ taking $L^2(\rho_1)$ to
$L^2(\rho_2)$ is given by
%
\begin{equation}
\label{eqrank-one}
\| u \otimes v \| =
\| u \|_{2, \rho_2} \| v \|_{2, \rho_1^{-1}}.
\end{equation}
Here the norm can be trace, Hilbert--Schmidt or operator norm, since all
agree for a rank-one operator. 

The finite rank terms include ones of the form $\tilde L( a \otimes
\tilde\varepsilon b)$. We use~\eqref{eqrank-one} to establish entrywise
bounds
%
\begin{equation}
\label{eqLbd1}
\pmatrix{
\| a \otimes\tilde\varepsilon b \| & \| a \otimes b \| \vspace*{2pt}\cr
\| \tilde\varepsilon a \otimes\tilde\varepsilon b \|
& \| \tilde\varepsilon b \otimes a \|
}
\leq
\pmatrix{
A_+ B_- & A_+ B_+ \vspace*{2pt}\cr A_- B_- & A_+ B_-},
\end{equation}
where
\begin{eqnarray*}
A_+ & = &\| a \|_+,\qquad B_+ = \| b \|_+, \\
A_- & =& \| \tilde\varepsilon a \|_-,\qquad
 B_-  = \| \tilde\varepsilon b \|_-.
\end{eqnarray*}
%
Indeed, for the
$(i,j)$th entry, apply~\eqref{eqrank-one} to $ a_{ij}
\otimes b_{ij}\dvtx  L^2(\rho_j) \to L^2(\rho_i)$, where $\rho_1 = \rho$
and $\rho_2 = \rho^{-1}$. On the right, and henceforth, we abbreviate
the $L^2$ norms on $L^2(\rho)$ and $L^2(\rho^{-1})$ by $\| \cdot\|_+$
and $\| \cdot\|_-$, respectively.



Let us indicate how this applies to
%
\[
- 2 \delta_0^F =
\tilde L [\phi_\tau\otimes\tilde\varepsilon(\psi_\tau-
A_N) + (\phi_\tau- A) \otimes\tilde\varepsilon A_N].
\]
Consider the first term on the right-hand side---the second term is
similar---and apply~\eqref{eqLbd1} with $a = \phi_\tau, b = \psi
_\tau
- A_N$. From Corollary~\ref{cororthogonal-case-1} we have
\begin{eqnarray*}
A_+^2 & =& \| \phi_\tau\|_+^2
= \int_{s_0}^\infty\phi_\tau^2 \rho
\leq C(\gamma) e^{-2s_0 + \gamma|s_0|}, \\
B_-^2 & = &\| \tilde\varepsilon(\psi_\tau- A_N) \|_-^2
\leq C(\gamma) N^{-4/3} e^{-s_0 + \gamma|s_0|}
\end{eqnarray*}
and with similar bounds, respectively, for $A_-^2$ and $B_+^2$.
Hence
%
\begin{equation} \label{eqdelF0}
A_\pm B_\pm\leq
C(\gamma) N^{-2/3} e^{-3s_0/2 + \gamma|s_0|}.
\end{equation}




Turning to the the $\delta_1^F, \delta_2^F$ terms, we have
%
\[
2 \delta_1^F =
\pmatrix{
(u_{N1} - A) \otimes1 & 0 \vspace*{2pt}\cr
-(u_{N2} - \tilde\varepsilon A) \otimes1 & 0}
,
\qquad
2 (\delta_2^F)^t =
\pmatrix{
0 & 1 \otimes(u_{N2} - \tilde\varepsilon A) \vspace*{2pt}\cr
0 & 1 \otimes(u_{N1} - A)}
\]
with $ u_{N1} = \beta_{N-1} \phi_\tau$ and
$ u_{N2} = \beta_{N-1} \tilde\varepsilon\phi_\tau$.
Using~\eqref{eqLbd1}, we find that the norms of the terms in the
first column of $\delta_1^F$ are bounded by
$ \| u_{N1} - A \|_+ \| 1 \|_-$ and
$ \| u_{N2} - \tilde\varepsilon A \|_- \| 1 \|_- $
while the norms of the second column of $(\delta_2^F)^t$ are bounded
by the same quantities interchanged.

From the definitions, and with $s_0 \geq0$, we have
$ \| 1 \|_-^2
\leq\gamma^{-1} e^{-\gamma s_0}$
and
\[
\| u_{N1} - A \|_+ \leq
|\beta_{N-1} - 1| \| \phi_\tau\|_+ +
\| \phi_\tau- A \|_+.
\]
Note that $|\beta_{N-1} - 1| = O(N^{-1})$.
Using also the bounds of Corollary~\ref{cororthogonal-case-1},
%
\[
\|u_{N1} - A\|_+
\leq(C N^{-1} e^{-s_0} + CN^{-2/3} e^{-s_0/2}
)e^{\gamma s_0/2}
\]
and $ \|u_{N1} - A\|_+ \| 1 \|_- \leq C N^{-2/3} e^{-s_0/2}$.
The term $\| u_{N2} - \tilde\varepsilon A \|_- $ is bounded analogously.

We finally assemble the bounds obtained from
(\ref{eqdelRI}), (\ref{eqdelF0}) and the analysis of $\delta_i^F$
and only track the tail dependence on $s_0$ for $s_0 > 0$.
Thus (\ref{eqF1summary}) is bounded by
%
\[
C N^{-2/3} ( e^{-3s_0/2 + \gamma s_0} + e^{-s_0/2}),
\]
where the second term results from $\delta_1^F$ and $\delta_2^F$.
It is clear that
$\gamma= 1$ yields a bound $C N^{-2/3} e^{- s_0/2}.$


\begin{appendix}\label{app}
\section{Hermite polynomial asymptotics near largest~zero}\label{secjacobi-asy}




Define new independent and dependent variables $\zeta$ and $W$ via the
equations~\eqref{eqnewvarsa},
which put~\eqref{eqhermite-de} into the form
%
\begin{equation}
\label{eqkappaform}
\frac{d^2 W}{d\zeta^2} = \{ \kappa^2 \zeta+ \psi(\zeta) \} W,
\end{equation}
where the perturbation term
$\psi(\xi) = \dot\zeta^{-1/2} ( d^2/ d \zeta^2 ) (\dot\zeta^{1/2}).$
If the perturbation term $\psi(\zeta)$ were
absent, the equation $d^2 W / d \zeta^2 = \kappa^2 \zeta W$ would have
linearly independent\vadjust{\goodbreak} solutions in terms of the Airy functions
$\operatorname{Ai}(\kappa^{2/3} \zeta)$ and
$\operatorname{Bi}(\kappa^{2/3} \zeta)$. Our interest is in
approximating the
\textit{recessive} solution $\operatorname{Ai}(\kappa^{2/3} \zeta
)$, so
write the relevant solution of (\ref{eqkappaform}) as $W_2(\zeta) =
\operatorname{Ai}( \kappa^{2/3} \zeta) + \eta(\zeta).$
In terms of the original independent and dependent variables $w$
and $\xi,$ the solution~$W_2$ becomes
%
\begin{equation}
w_2(\xi, \kappa) = \dot\zeta^{-1/2}(\xi) \{ \Ai(\kappa
^{2/3}\zeta) +
\varepsilon_2 (\xi, \kappa) \}.
\label{eqw2expr}
\end{equation}

Olver (\citeyear{olve74})---hereafter abbreviated as [O]---provides,
in his Theorem 11.3.1,
an explicit bound for $\eta(\zeta)$ and hence $\varepsilon_2$
and its derivative.
To describe these error bounds even in the oscillatory region of $\Ai
(x)$, [O] introduces a~positive weight function $\mathsf{E}(x)\geq1$ and
positive moduli functions $\mathsf{M}(x) \leq1$ and $\mathsf{N}(x)$
such that
for all $x$,
%
\begin{equation}
\label{eqAiAiprimebd}
|\Ai(x)| \leq\mathsf{M}(x) \mathsf{E}^{-1}(x),\qquad
|\Ai'(x)| \leq\mathsf{N}(x) \mathsf{E}^{-1}(x).
\end{equation}
[Here, $\mathsf{E}^{-1}(x)$ denotes $1/ \mathsf{E}(x).$]
In addition,
%
\begin{equation}
\Ai(x) = 2^{-1/2} \mathsf{M}(x) \mathsf{E}^{-1}(x),\qquad
x \geq c \doteq-0.37,
\label{eqmodweight}
\end{equation}
and the asymptotics as $x \rightarrow\infty$ are given by
%
\begin{equation}
\qquad\mathsf{E}(x) \sim\sqrt{2} e^{({2}/{3})x^{3/2}}, \qquad\mathsf{M}(x)
\sim
\pi^{-1/2} x^{-1/4}\quad \mbox{and}\quad
\mathsf{N}(x) \sim\pi^{-1/2} x^{1/4}.
\label{eqemasymp}
\end{equation}
The key bounds of [O, Theorem 11.3.1] then state, for $\xi> 0$ and
$\hat f(\xi) = f(\xi)/\xi$,
%
\begin{eqnarray}
|\varepsilon_2 (\xi, \kappa)| & \leq&
(\mathsf{M}/\mathsf{E})(\kappa^{2/3}\zeta)
\biggl[\exp\biggl\{\frac{\lambda_0}{\kappa}
\mathcal{V}(\zeta)\biggr\} - 1\biggr], 
\label{eqkeybound} \\
| \partial_\xi\varepsilon_2 (\xi, \kappa)| & \leq&
\kappa^{2/3} N^{-1} \hat f^{1/2}(\xi)
(\mathsf{N}/\mathsf{E})(\kappa^{2/3}\zeta),
\label{eqkeyboundprime}
\end{eqnarray}
where $\lambda_0 \doteq1.04$.
For $\kappa^{2/3} \zeta\geq c$,
\eqref{eqmodweight} shows that the coefficient
in~\eqref{eqkeybound} is just $\sqrt{2} \Ai( \kappa^{2/3}\zeta)$.
Here $ \mathcal{V}(\zeta) = \mathcal{V}_{[\xi,\infty]}(H)$
is the total variation on $[\xi,\infty]$ of the \textit{error control}
function
$H(\xi) = - \int_0^{\zeta(\xi)} |v|^{-1/2} \psi(v) \,dv.$
From [O, page~403] we have
$\lambda_0 \mathcal{V}_{[\xi,\infty)}(H) \leq0.28$
and hence
%
\begin{equation}
\label{eqNm1bd}
\exp\biggl\{ \frac{\lambda_0}{\kappa} \mathcal{V}(\zeta) \biggr\} - 1 \leq1/N.
\end{equation}

\textit{Application to Hermite polynomials}.
In the case of Hermite polynomials,
transformed as in~\eqref{eqhermite-de}, the points
$ \pm\infty$ are irregular singularities, and the points $\xi_\pm=
\pm1$ are turning points.
We are interested in behavior near the upper turning point $\xi_+$,
which is located near the largest (scaled) zero of~$H_N$.
Using (\ref{eqnewvarsb}), the independent variable
$\zeta(\xi)$ is given in terms of $f(\xi)$ by
%
\begin{equation}
\label{eqzetadef}
(2/3) \zeta^{3/2}(\xi)
 = \tfrac{1}{2} \xi(\xi^2 -1)^{1/2} - \tfrac{1}{2} \log\bigl( \xi+ (\xi^2
-1)^{1/2}\bigr)
\end{equation}
for $\xi\geq1$, and by
\[
(2/3) (-\zeta)^{3/2}(\xi)
 = \tfrac{1}{2} [\cos^{-1} \xi- \xi(1 - \xi^2)^{1/2}],
\]
for $\xi\leq1.$
The function $\zeta(\xi)$ is increasing and $C^2$ on $(0,\infty)$
(e.g., [O, page~399]), with $\ddot\zeta(\xi)$ nonnegative and bounded.
It is easily seen that $\zeta\rightarrow\infty$
as $\xi\rightarrow\infty,$
and more precisely, from~\eqref{eqzetadef},
that
%
\begin{equation}
\label{eqxitoinf}
(2/3) \zeta^{3/2}(\xi) = \tfrac{1}{2}\bigl(\xi^2 - \log\xi- \tfrac
{1}{2} -
\log2 \bigr) + O(\xi^{-2}),
\end{equation}
from which it follows that
%
\begin{equation}
\label{eqxidotlim}
\dot\zeta(\xi) \sim(4 \xi/3)^{1/3}\qquad
\mbox{as } \xi\rightarrow\infty.
\end{equation}

We remark that $\dot\zeta= \dot\zeta(1)$ is easily evaluated using
L'H\^opital's rule. From (\ref{eqzetadef}), as $\xi\rightarrow1$,
we have
$ \dot\zeta^2 (\xi) = (\xi^2 - 1)/\zeta(\xi)
\rightarrow2/\dot\zeta,$
so that $\dot\zeta= 2^{1/3}.$ In addition, we shall need the
function
\[
r(\xi) = [ \dot\zeta(\xi) / \dot\zeta]^{-1/2},
\]
which is positive on $(0,\infty)$ since $\zeta(\xi)$ is strictly increasing.
Both $r(\xi)$ and~$\dot r(\xi)$ are continuous on $[0,\infty)$, and as
$\xi\rightarrow\infty$ we have
$r(\xi) \sim(2\xi/3)^{-1/6}$ and
$\dot r(\xi) \sim c_1 \xi^{-7/6}$,
so that $r(\xi)$ and $r'(\xi)$ are both bounded on $[0,\infty)$.

Bound (\ref{eqkeybound}) has a double asymptotic property in $\xi$ and
$\kappa$ which will be useful.
First, suppose that $N$, and hence $\kappa$, are held fixed.
As $\xi\rightarrow\infty, \mathcal{V}(\zeta) \rightarrow0$ and so from
(\ref{eqkeybound}) and its following remarks $\varepsilon_2(\xi
,\kappa) =
o( \A(\kappa^{2/3} \zeta)).$ Consequently, as $\xi\rightarrow
\infty$
we have
$ w_2(\xi,\kappa) \sim\dot\zeta^{-1/2} (\xi) \A(\kappa^{2/3}
\zeta).$
If the weighted polynomial $w_N(\xi)$ is a recessive solution of
\eqref{eqhermite-de}, then it must be proportional to $w_2$,
so that $ w_N(\xi) = c_N w_2(\xi,\kappa)$.
Now $c_N$ may be identified by comparing the growth of $w_N(\xi)$
as $\xi\rightarrow\infty$ with that of $w_2(\xi,\kappa)$
(Appendix~\ref{secidentification-cn}),
%
\begin{equation}
\label{eqcNbd}
c_N = e^{\theta''/N} \kappa_N^{1/6} (2/N)^{1/4},
\end{equation}
where $\theta'' = O(1)$.
Now we can use (\ref{eqw2expr}) to write
$\phi_N(x) = w_N(\xi)$ in terms of the Airy approximation.
Below, we write $\bar e_N$ for any term, that is,
uniformly $1 + O(N^{-1})$. Hence
\[
(2N)^{1/4} \phi_N(x)
= \bar e 2^{1/2} \kappa_N^{1/6} w_2(\xi,\kappa).
\]
Set $N_+ = N+1/2$ and $\bar\tau_N = 2^{-1/2} N_+^{-1/6}$.
Since $2^{1/2} \kappa_N^{1/6} \dot\zeta(\xi)^{-1/2} = \bar
\tau_N^{-1} r(\xi)$, and
using the Airy approximation~\eqref{eqw2expr} to $w_2(\xi,\kappa)$,
we finally have
%
\begin{equation}
\label{eqbasic-rep}
(2N)^{1/4} \bar\tau_N \phi_N(x) = \bar e_N r(\xi)
\{ \Ai(\kappa_N^{2/3} \zeta) + \varepsilon_2(\xi,\kappa_N) \}.
\end{equation}
%

\textit{Approximations at degree $N$ and $N-1$}.
The kernel $S_{N,2}(x,y)$ is expressed
in terms of the two functions $\phi_{N-1}(x)$ and $\phi_{N}(x)$,
which need separate Liouville--Green asymptotic approximations.
Thus, for example, in comparing the two cases, we have $\kappa_N =
2N+1$ and $\kappa_{N-1} = 2N-1$. The turning point $\xi_+ = 1$ and the
transformation $\zeta(\xi)$ of (\ref{eqzetadef}) are the same in both
cases, hence so is $r(\xi)$. The analog of (\ref{eqbasic-rep}) then
states
%
\[
(2N-2)^{1/4} \bar\tau_{N-1} \phi_{N-1}(x) = \bar e_{N-1} r(\xi)
\{ \Ai(\kappa_{N-1}^{2/3} \zeta) + \varepsilon_2(\xi,\kappa
_{N-1}) \}.
\]

Rather than $\bar\tau_{Nj} = 2^{-1/2} N_\pm^{-1/6}$, we will use the
single factor $\tau_N = 2^{-1/2} N^{-1/6}$ in the work below.
Clearly, we may replace both $(2N)^{1/4} \bar\tau_N$
in~\eqref{eqbasic-rep} and
$(2N-2)^{1/4} \bar\tau_{N-1}$ in the preceding display by $(2N)^{1/4}
\tau_{N}$ at cost of multiplicative error terms
$ e_{Nj} = 1 + O(N^{-1})$.

To summarize then, with the convention that quantities with
subscript~$Nj$ differ for $Nj = N, N-1$, while those with subscript $N$ do not,
we have
%
\begin{equation}
\label{eqbasic-rep-mod}
(2N)^{1/4} \tau_N \phi_{Nj}(x) = e_{Nj} r(\xi)
\{ \Ai(\kappa_{Nj}^{2/3} \zeta) + \varepsilon_2(\xi,\kappa_{Nj})
\}.
\end{equation}
%

Denote the left-hand side of~\eqref{eqbasic-rep-mod} by $\bar\phi_{Nj}$.
We seek a uniform bound on the Airy approximation.
If we write $x = \sqrt{\kappa_{Nj}}\xi$ in the form $u_{Nj} + s
\tau_N$, then we have in particular
$u_N = \sqrt{2N+1}$ and $u_{N-1} = \sqrt{2N-1}$.
In turn,
%
\[
\xi= 1 + s \tau_N/ \sqrt{\kappa_{Nj}} = 1 + s \sigma_{Nj},
\]
where we define
%
\begin{equation}
\label{eqsigNdef}
\qquad\sigma_{Nj} = \tau_N/u_{Nj}
= 2^{-1/2} N^{-1/6} (2N \pm1)^{-1/2}
= 2^{-1} N^{-2/3} \bigl(1+O(N^{-1})\bigr).
\end{equation}

We turn now to the proof of Proposition~\ref{prop-localbds}.
We first record some properties of the map
$s \rightarrow\kappa_{Nj}^{2/3} \zeta(1 + \sigma_{Nj}s)$, which we
sometimes abbreviate as $\kappa^{2/3} \zeta$.

\begin{lemma}
\label{lemlocalexp}
Given $s_L \in\mathbb{R}$,
%
\begin{eqnarray}
|\kappa^{2/3} \zeta- s| & \leq& |s|/4\qquad
\mbox{for }  s_L \leq s \leq N^{1/6}, N \geq N_0,
\label{eqonesixth} \\
| \kappa^{2/3} \zeta| & \leq& C |s|/4 \qquad
\mbox{for }  s_L \leq s \leq C N^{2/3}, \mbox{ all } N.
\label{eqthreehalves}
\end{eqnarray}
\end{lemma}
\begin{pf}
Expand $\zeta(\xi)$ about the turning point $\xi_+ = 1$:
%
\begin{equation}
\label{eqtaylor}
\kappa_{Nj}^{2/3} \zeta(1 + \sigma_{Nj} s) =
\kappa_{Nj}^{2/3} \sigma_{Nj} \dot\zeta s +
\tfrac{1}{2}\kappa_{Nj}^{2/3} \sigma_{Nj}^2 s^2 \ddot\zeta(\xi^*).
\end{equation}
We note from the definitions that
%
\begin{equation}
\label{eqoneplus}
\kappa_{Nj}^{2/3} \sigma_{Nj} \dot\zeta
= \bigl(1 \pm1/(2N)\bigr)^{1/6} = 1 + \delta_N,
\end{equation}
with $ | \delta_N| \leq N^{-1}$ for all $N \geq1.$
Since $0 \leq\ddot\zeta$ is bounded,
we find that
%
\begin{equation}
| \kappa^{2/3} \zeta- s|
\leq\biggl(\frac{1}{N} + \frac{M |s|}{N^{2/3}}\biggr)
|s|, \label{eqlocalbdd}
\end{equation}
%
again for all $N \geq1.$
If $s < N^{1/6}$, then the right-hand side is bounded by $|s|/4$ for $N
\geq
N_0(M,s_L)$.
If $|s| < N^{2/3}$, then we have (\ref{eqthreehalves}) for $C =
C(s_L,M)$.
\end{pf}

We consider some global bounds, valid for $s \geq s_L$, or
equivalently for $\xi\geq1 + s_L \sigma_{Nj}$.

\begin{lemma}
\label{propE-inv-bound} Let $s_L < 0.$
Let $\xi= 1 + \sigma_{Nj} s$ with $\sigma_{Nj}$ satisfying
\eqref{eqsigNdef}.
There exists $C = C(s_L)$ such that for $s \geq s_L$,
\begin{eqnarray*}
\mathsf{E}^{-1} (\kappa_{Nj}^{2/3} \zeta)
& \leq& C e^{-2s}, 
\\
\mathsf{N}(\kappa^{2/3} \zeta)
& \leq &C( 1 + |s|^{1/3}). 
\end{eqnarray*}
\end{lemma}

Some immediate consequences: using~\eqref{eqAiAiprimebd}
and $\mathsf{M} \leq1$, for $s \geq s_L$,
%
\begin{eqnarray}
| A(\kappa_{Nj}^{2/3} \zeta(\xi))|
& \leq&| \mathsf{M}/\mathsf{E} |(\kappa^{2/3} \zeta)
\leq C e^{-2s}, \label{eqAbound} \\
\label{eqcoroll1}
|\varepsilon_2(\xi, \kappa)|
& \leq& N^{-1} |\mathsf{N}/\mathsf{E}|(\kappa^{2/3} \zeta) \leq C
N^{-1} e^{-2s}, \\
| A'(\kappa_{Nj}^{2/3} \zeta(\xi))|
& \leq&| \mathsf{N}/\mathsf{E} |(\kappa^{2/3} \zeta)
\leq C (1 + |s|^{1/3}) e^{-2s} \label{eqcoroll2}.
\end{eqnarray}

\begin{pf}
First, since
$f(\xi) = (\xi+ 1)(\xi- 1) \geq2 \sigma_{Nj} s$,
we use (\ref{eqoneplus})
to observe that for $s \geq r^2$,
%
\[
\kappa_{Nj} \sigma_{Nj} \sqrt{f}
\geq\sqrt{2} \kappa_{Nj} \sigma_{Nj}^{3/2} \sqrt{s}
\geq e_{Nj} r.
\]
Hence, from (\ref{eqnewvarsb}), again for $s \geq r^2$,
%
\[
\frac{2}{3} \kappa_{Nj} \zeta^{3/2}
= \kappa_{Nj} \int_{1}^\xi\sqrt{f}
= e_{Nj} r (s - r^2).
\]
Now choose $r$ large enough so that for $N>N_0$ and $j = N, N-1$, we
have \mbox{$e_{Nj}r \geq1$}.
From~\eqref{eqemasymp} we have $\mathsf{E}^{-1}(s) \leq C \exp( -
\frac{2}{3}
s^{3/2} )$ for $s \geq0,$ and so in particular for $s \geq r^2$,
%
\[
\mathsf{E}^{-1}(\kappa_{Nj}^{2/3} \zeta) \leq C(r) e^{-s}.
\]
For $s \in[s_L,r^2],$ we simply use the bound $\mathsf{E} \geq1.$

For the second statement,
we will use the bound $\mathsf{N}(s) \leq1 + |s|^{1/4}$ [O, pages~396--397].
First, for $s \leq N^{2/3}$, using the bound on $\mathsf{N}$ and
(\ref{eqthreehalves}), we obtain
$ \mathsf{N}(\kappa^{2/3} \zeta) \leq C( 1 + |s|^{1/4})$.
When $s \geq N^{2/3}$, we use (\ref{eqzetadef}) to bound
%
\[
(2/3) \zeta^{3/2} (\xi) \leq\int_1^{\xi} t \,dt \leq\xi^2/2
\leq1 + \sigma_{Nj}^2 s^2
\leq c_0 \sigma_{Nj}^2 s^2.
\]
From (\ref{eqsigNdef}) we have $\kappa_{Nj} \sigma_{Nj}^2 = 2
N^{-1/3}$ and so $\kappa_{Nj} \zeta^{3/2} \leq c_1 s^2$ for all $N$
and hence $\mathsf{N}(\kappa^{2/3} \zeta) \leq1 + c_1^{1/6} s^{1/3}$
as required.
\end{pf}



\begin{pf*}{Proof of Proposition \protect\ref{prop-localbds}}
We begin from the formula
%
\begin{equation}
\label{eqinitial}
\bar\phi_{Nj}(u_{Nj} + s \tau_N)
= e_{Nj} r(\xi) \{ \Ai(\kappa_{Nj}^{2/3} \zeta) + \varepsilon
_2(\xi,
\kappa_{Nj}) \}.
\end{equation}
The bound (\ref{eqphiNjbd}) then follows from~\eqref{eqcoroll1},
\eqref{eqAbound} and boundedness of $r(\xi)$.
To ease notation, we will, as needed, drop subscripts from $e_{Nj},
\sigma_{Nj},
\kappa_{Nj},$ and~$\tau_N$, writing~$\bar e$ for a term, that is,
generically $1 + O(N^{-1})$.\vadjust{\goodbreak}

For the next bound, we differentiate (\ref{eqinitial}), obtaining
%
\begin{equation}
\label{eqphiNdiff}
\tau_N \bar\phi_{Nj}^\prime(u_{Nj} + s \tau_N)
= D_1 + D_2 + D_3,
\end{equation}
with the component terms given by
\begin{eqnarray*}
D_1 & =& \bar e \sigma\dot r(\xi) [\Ai(\kappa^{2/3} \zeta) +
\varepsilon
_2(\xi,
\kappa)], \\
D_2 & = &\bar e r(\xi) \Ai'( \kappa^{2/3} \zeta) \sigma\kappa^{2/3}
\dot\zeta(\xi), \\
D_3 & =& \bar e r(\xi) \sigma\partial_\xi\varepsilon_2(\xi, \kappa).
\end{eqnarray*}

Since $\dot r(\xi)$ is bounded, we again use (\ref{eqcoroll1}) and
\eqref{eqAbound} to conclude that
$|D_1| \leq C \sigma_N e^{-2s}$ for all $s \geq s_L$.
From~\eqref{eqoneplus}
and (\ref{eqxidotlim}), we observe that
%
\begin{equation}
\label{eqxidbd}
\sigma_N \kappa_N^{2/3} \dot\zeta(\xi) \leq C |\xi|^{1/3}
\leq C (1 + \sigma_N^{1/3} |s|^{1/3}).
\end{equation}

Turning to the second term, we have from (\ref{eqxidbd}) and
\eqref{eqcoroll2} that
%
\[
|D_2| \leq
C |\Ai'(\kappa^{2/3} \zeta) \sigma\kappa^{2/3} \dot\zeta(\xi)|
\leq C (1 + \sigma_N^{1/3} |s|^{1/3}) (1 + |s|^{1/3}) e^{-2s}
\leq C e^{-s}.
\]

Using~\eqref{eqkeyboundprime} and~\eqref{eqNm1bd}, we can rewrite
$D_3$ as
%
\[
|D_3| \leq
C N^{-1} \cdot r(\xi) \sigma_N \kappa_N^{2/3} \dot\zeta(\xi)
\cdot(\mathsf{N}/\mathsf{E})(\kappa^{2/3} \zeta).
\]
Using also boundedness of $r(\xi)$,~\eqref{eqxidbd}
and~\eqref{eqcoroll2},
we conclude for all $s \geq s_L$,
%
\[
|D_3| \leq C N^{-1} \cdot\sqrt2 \cdot(1 + \sigma_N^{1/3}
|s|^{1/3}) ( 1 + |s|^{1/3}) e^{-2s} \leq C N^{-2/3} e^{-s}.
\]
This completes the proof of bound (\ref{eqphiNjbd}) for $\tau_N \bar
\phi_{Nj}$.
\end{pf*}

Turning to the error bound (\ref{eqdelphiNj}) and its analog for $\tau
_N \bar
\phi_{Nj}$,
we first note that we may confine attention to
$s \in[s_L, N^{1/6}]$, since for $s \geq N^{1/6}$, the bounds follow
trivially from (\ref{eqphiNjbd}) and its analog
and~\eqref{eqAiFacts}.

We use the decomposition suggested by~\eqref{eqinitial},
\begin{eqnarray*}
\bar\phi_{Nj}(x) - \A(s)
& = &[ e_{Nj} r(\xi) - 1] \Ai( \kappa_{Nj}^{2/3} \zeta)
+ [\Ai(\kappa_{Nj}^{2/3} \zeta) - \Ai(s)] \\
&&{}+ e_{Nj} r(\xi) \varepsilon_2(\xi,\kappa_{Nj}) \\
& =& E_{N1} + E_{N2} + E_{N3}.
\end{eqnarray*}

For the $E_{N1}$ term, first use
$\xi= 1 + \sigma_{Nj} s$ to write
%
\begin{equation}
\label{eqrat-minus-1}
| \dot\zeta(\xi)/ \dot\zeta- 1 |
= \biggl| \int_1^{1 + s \sigma_{Nj}} \ddot\zeta(u) / \dot\zeta \,du \biggr|
\leq C s \sigma_{Nj},
\end{equation}
since $\ddot\zeta(u)$ is bounded for $u \in[1, 1 + s \sigma_{Nj}]
\subset
[1, 1+ N^{-1/2}]$.
Together with $r(\xi) = [\dot\zeta(\xi)/ \dot\zeta]^{-1/2}$, this yields
%
\begin{equation}
\label{eqeNrN-1}
| e_{Nj} r(\xi) - 1| \leq C(1+ s) \sigma_{Nj}.
\end{equation}
Combined with~\eqref{eqAbound}, we obtain
$ |E_{N1}| \leq C \sigma_{Nj} (1+ s) e^{-2s} \leq C N^{-2/3}
e^{-s/2}.$\vadjust{\goodbreak}

For the $E_{N2}$ term, we use (\ref{eqonesixth}) and
(\ref{eqlocalbdd}) to write
\begin{eqnarray*}
| \Ai(\kappa_{Nj}^{2/3} \zeta) - \Ai(s)|
& \leq& C |s| N^{-2/3} (N^{-1/3} + |s|) \sup\bigl\{ | \Ai'(t)| \dvtx  \tfrac
{3}{4} s
\leq t \leq\tfrac{5}{4}s \bigr\} \\
& \leq& C N^{-2/3} e^{-s/2},
\end{eqnarray*}
uniformly for $s \in[s_L, N^{1/6}]$,
where we used~\eqref{eqAiFacts}.

Finally, for the $E_{N3}$ term,~\eqref{eqcoroll1}
and boundedness of $r$ imply that
\[
|E_{N3}| \leq e_{Nj} r(\xi) (\mathsf{M}/\mathsf{E}) (
\kappa_{Nj}^{2/3} \zeta) N^{-1} \leq C N^{-1} e^{-2s}.
\]

We turn now to the proof of (\ref{eqdelphiNj}) for $\tau_N \bar
\phi_{Nj}$ on $[s_L, N^{1/6}]$.
Using~\eqref{eqphiNdiff}, we may write the difference
$ \tau_N \bar\phi_{Nj}'(x) - \Ai'(s) $ as
$D_1(s) + [D_2(s) - \A'(s)] + D_3(s)$.
Now decompose $D_2(s) - \A'(s) = G_1 + G_2 + G_3$, with
%
\[
G_{1} = [\bar e r(\xi) - 1][ \dot\zeta(\xi)/ \dot\zeta]
\Ai'(\kappa^{2/3} \zeta),\qquad
G_{2} = [ \dot\zeta(\xi)/ \dot\zeta-1 ]
\Ai'(\kappa^{2/3} \zeta)
\]
and $G_{3} = \Ai'(\kappa^{2/3} \zeta) - \Ai'(s)$.
Combining~\eqref{eqeNrN-1},~\eqref{eqcoroll2} and then
\eqref{eqrat-minus-1}, we find
\begin{eqnarray*}
|G_{1}| & \leq& C(1+s) \sigma_N \cdot2 \cdot C ( 1 + |s|^{1/3}) e^{-2s}
\leq C N^{-2/3} e^{-s/2}, \\
|G_{2}| & \leq& C s \sigma_N \cdot C ( 1 + |s|^{1/3}) e^{-2s}
\leq C N^{-2/3} e^{-s/2}.
\end{eqnarray*}
$G_{3}$ is treated in exactly the same manner as the $E_{N2}$ term
above, additionally using the equation $\Ai''(x) = x \Ai(x)$.

\section{Identification of $c_N$}
\label{secidentification-cn}

We first remark that as $\zeta\rightarrow\infty$ when $\xi
\rightarrow
\infty$, we may substitute
$\Ai(x) \sim[2 \sqrt\pi x^{1/4}]^{-1} \exp\{ -(2/3) x^{3/2} \}$
into (\ref{eqw2expr}), along with $\dot\zeta^{-1/2} =
[\zeta/f(\xi)]^{1/4}$ from (\ref{eqnewvarsa}) to obtain
\[
w_2(\xi,\kappa) \sim
[2 \sqrt\pi]^{-1} \kappa^{-1/6} f^{-1/4}(\xi)
\exp\{ -(2/3) \kappa\zeta^{3/2} \}.
\]
Consequently
we may express $c_N$ in terms of the limit
%
\[
c_N = \lim_{\xi\rightarrow\infty} w_N(\xi) \cdot2 \sqrt\pi
\kappa^{1/6} f^{1/4}(\xi) \cdot\exp\{ (2/3) \kappa\zeta^{3/2} \}.
\]

Write $N_+$ for $N+1/2.$
Since $w_N(\xi) = \phi_N(x)$ with
$x = \sqrt{2 N_+} \xi$, and since
$H_N(x) \sim2^N x^N$, we have as $\xi\rightarrow\infty$,
%
\[
w_N(\xi) = h_N^{-1/2} e^{-N_+ \xi^2} H_N\bigl( \sqrt{2N_+} \xi\bigr)
\sim h_N^{-1/2} e^{-N_+ \xi^2 + N \log\xi} 2^N (2N_+)^{N/2},
\]
and $ f^{1/4}(\xi) = e^{ (\log\xi)/2 + O(\xi^{-2})}$,
while from~\eqref{eqxitoinf}
%
\[
\exp\{ (2/3) \kappa\zeta^{3/2} \}
= e^{N_+ \xi^2 - N_+ \log\xi- N_+/2 - N_+ \log2 + O(\xi^{-2})}.
\]

Multiply the last three quantities: the coefficients of $\xi^2$ and
$\log\xi$ cancel, leaving $\xi$-dependence of only $O(\xi^{-2})$ as
$\xi
\rightarrow\infty$.
Hence
%
\[
c_N = 2 \sqrt{\pi} \kappa^{1/6} h_N^{-1/2} (2N_+)^{N/2} e^{-N_+/2} 2^{N-N_+},
\]
and noting that $(N/2) \log N_+ = (N/2) \log N + 1/4 + O(N^{-1})$, we get
%
\[
\kappa^{-1/6} c_N \sqrt{h_N}
= \sqrt{2 \pi} \exp\biggl\{ \frac{N}{2} \log2 + \frac{N}{2} \log N -
\frac{N}{2} + O\biggl(\frac{1}{N}\biggr) \biggr\}.
\]
%
Applying Stirling's formula to $h_N = \sqrt{\pi} 2^N N!$,
and dividing into the previous display yields (\ref{eqcNbd}).


\end{appendix}

\section*{Acknowledgments}
We thank Leonard Choup, Peter Forrester, Craig Tracy
and Nicholas
Witte for very helpful discussions and Folkmar Bornemann for providing
software. Iain M. Johnstone
thanks University of Barcelona for hospitality during writing of the paper.


%

\printaddresses

\end{document}